\documentclass[12pt,a4paper]{article}
\usepackage[T2A]{fontenc}
\usepackage[cp1251]{inputenc}
\usepackage[russian, english]{babel}
\usepackage{epsfig,amssymb,color,graphics}
\usepackage{amsthm, amsmath, amssymb, amsbsy, amsfonts, latexsym, euscript}
cm \hoffset = -1.0 true cm

\newtheorem{Th}{Theorem}
\newtheorem{st}{Statement}[section]
\newtheorem{Cor}{Corollary}[section]
\newtheorem{Prop}{Proposition}[section]
\newtheorem{Lemm}{Lemma}[section]
\newtheorem{Defin}{Definition}[section]

\newenvironment{demo}{{\bf Proof: }}
{\hfill $\diamond$\medskip}

\begin{document}

\title{
\medskip
\Large On Embedding of Multidimensional Morse-Smale Diffeomorphisms into Topological Flows}
\author{V.~Grines, E.~Gurevich, O.~Pochinka}
\date{}

\maketitle
\abstract{J.~Palis found necessary  conditions for  a Morse-Smale diffeomorphism on a closed $n$-dimensional manifold $M^n$  to  embed into a topological  flow and proved that these conditions are also sufficient for $n=2$.  For the case $n=3$ a possibility of  wild embedding of closures of separatrices of saddles is an additional obstacle  for Morse-Smale cascades to embed into topological flows.  In this paper we show that there are no such obstructions   for  Morse-Smale diffeomorphisms without heteroclinic intersection given  on  the sphere $S^n, \,n\geq 4$,  and Palis's conditions again are  sufficient for such  diffeomorphisms.} 

\sloppy

\section{Introduction and statements of results}

\hspace*{\parindent} 

Let $M^n$ be  a smooth connected  closed $n$-manifold. Recall that  {\it  a $C^m$-flow} ($m\geq 0$) on the manifold  $M^n$  is a continuously depending on  $t\in\mathbb R$ family of  $C^m$-diffeomorphisms  $X^t:M^n\to M^n$ that satisfies the following conditions: 
\begin{itemize}
\item[1)] $X^0(x)=x$ for any point  $x\in M^n$;
\item[2)] $X^t(X^s(x))=X^{t+s}(x)$ for any $s,t\in \mathbb{R}$, $x\in M^n$. 
\end{itemize}
A $C^0$-flow is also called a {\it topological flow}. 
One says that a homeomorphism (diffeomorphism)  $f:M^n\to M^n$ {\it embeds} into a $C^m$-flow on $M^n$ if $f$ is the time one map of this flow.

 Obviously, if a homeomorphism embeds in a flow then it is isotopic to identity. For a homeomorphism of the  line and a connected subset of the line this condition also is necessary  (see~\cite{FSch},\cite{F}). If an orientation preserving homeomorphism $f$ of the circle satisfies either one of the three conditions: 1) $f$  has  a fixed point, 2) $f$  has  a dense orbit, or 3) $f$   is periodic  then it embeds in a flow (see~\cite{FUtz}). Sufficient conditions of embedding  in topological flow for a  homeomorphisms of a compact two-dimensional disk and of the  plane  one can find in review~\cite{Utz81}. An analytical, $\varepsilon-$closed to the identity diffeomorphism  $f:M^n\to M^n$ can be approximated with accuracy  $e^{-\frac{c}{\varepsilon}}$  by a diffeomorphism which embeds in an analytical flow, see~\cite{Tr}.  

Due to~\cite{Pa74}  the  set of  $C^r$-diffeomorphisms  ($r\geq 1$) which embed in $C^1$-flows is a  subset of the first category in  $Diff^r(M^n)$. As  Morse-Smale diffeomorphisms are  structurally stable  (see \cite{Pa},~\cite{PaSm69})   then for  any manifold   $M^n$ there exists an open  set (in $Diff^1(M^n)$)  of Morse-Smale diffeomorphisms  embeddable  in topological flows. This  set contains  neighborhoods of time one maps of  Morse-Smale flows without periodic trajectories (according to~\cite{Sm} such flows  exist on an  arbitrary smooth manifold). 

 Recall that a diffeomorphism  $f:M^n\to M^n$ is called  {\it a  Morse-Smale diffeomorphism} if it satisfies the following conditions:
\begin{itemize}
\item the non-wandering set $\Omega_f$ is finite and consists of hyperbolic periodic points;
\item for any two points  $p,q\in \Omega_f$ the intersection of the stable manifold  $W^s_p$  of the point  $p$ and the unstable manifold   $W^u_q$ of the point   $q$ is transversal\footnote{Definitions of stable and unstable manifolds and of transversality are given in the section~\ref{MShom}; see also  the book~\cite{GrPo-book} for references.}.
\end{itemize}

In~\cite{Pa}  J.~Palis  established the following necessary conditions of the embedding of a  Morse-Smale  diffeomorphism $f:M^n\to M^n$ into a topological flow (we  call  them  {\it Palis conditions}):
\begin{itemize}
\item[(1)] the non-wandering set  $\Omega_f$  coincides  with the set of fixed points of $f$;
\item[(2)] the restriction of the diffeomorphism  $f$ to each invariant manifold of a fixed point   $p\in \Omega_f$  preserves the orientation of the manifold;
\item[(3)] if for two distinct saddle points  $p,q\in \Omega_f$ the intersection     $W^s_p\cap W^u_q$ is not empty then it  contains no compact connected components.   
\end{itemize}

\begin{figure}
\begin{center}\includegraphics[width=0.75\textwidth]{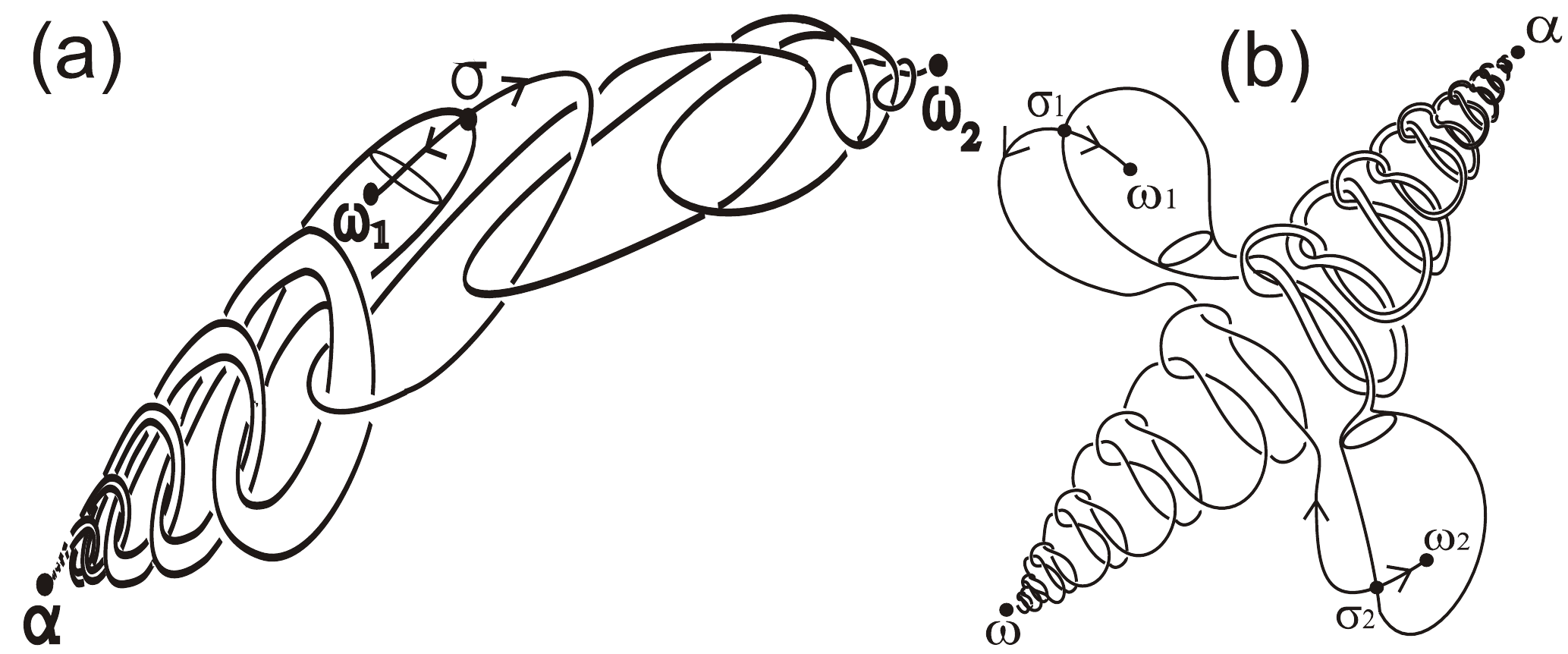}
\caption{Phase portraits of Morse-Smale diffeomorphisms on $S^3$ which do not   embed 
 in topological flows} \label{wi1}\end{center}
\end{figure}

According to  \cite{Pa}  these conditions are not only necessary but also sufficient for the case  $n=2$. For the case $n=3$ a possibility of  wild embedding of closures of separatrices of saddles is another  obstruction for Morse-Smale cascades to embed in topological flows (phase portraits of such diffeomorphisms are shown on the Figure~\ref{wi1}). In~\cite{GrGuMePo}  examples of such cascades are described  and  a   criteria  for embedding of Morse-Smale 3-diffeomorphisms in topological flows is provided.    In the present paper we  establish that  the Palis conditions are  sufficient  for  Morse-Smale diffeomorphisms on $S^n,\,n\geq 4,$ such that for any distinct saddle  points $p, q\in \Omega_f$  the intersection     $W^s_p\cap W^u_q$ is  empty.

\begin{Th}\label{emb-to-flow} Suppose that  a Morse-Smale diffeomorphism  $f: S^n\to S^n$, $n\geq 4$ satisfies the following conditions:
\begin{itemize}
\item[i)] the non-wandering set  $\Omega_f$  of the diffeomorphism  $f$ coincides with the set of its fixed points;
\item[ii)] the restriction of $f$ to each invariant manifold of a fixed point   $p\in \Omega_f$  preserves the orientation of the manifold;
\item[iii)] the invariant manifolds of distinct saddle points of $f$  do not intersect.

Then  $f$  embeds into a topological flow.
\end{itemize}
\end{Th}

\medskip
{\bf Acknowledgments.} Research is done with financial support of Russian Science Foundation  ({project  17-11-01041}) apart the section~\ref{aux}, which is done in frame of the  Basic Research  Program of HSE in  2018.

\section{Comments to  Theorem~\ref{emb-to-flow}}

Due to \cite{Pa} the conditions $i)$ and $ii)$ are necessary for embedding a Morse-Smale diffeomorphism into a flow. Our condition that the  ambient manifold  is the sphere $S^n$ and  the absence of heteroclinic intersections (condition iii)) are not necessary but  violation of each of them allows to construct  examples of Morse-Smale diffeomorphisms which do not embed in topological flows.  Below  we describe such examples.  

 In \cite{MeZh} V. Medvedev and E. Zhuzhoma constructed a  Morse-Smale diffeomorphism $f_0:M^4\to M^4$ satisfying conditions  $i)-iii)$ on a projective-like manifold $M^4$ (different from  $S^4$) whose  non-wandering set  consists of exactly three fixed points: a source, a sink and a saddle.  Invariant manifolds of the saddle are two-dimensional  and  the closure  of each of them  is a wild sphere (see \cite{MeZh},  Theorem~4, item~2). Assume that  $f_0$ embeds in a topological flow $X^t_0$. Then $X^t_0$ is a topological flow whose the non-wandering set consists of three equilibrium points with locally hyperbolic behavior. According to \cite[Theorem 3]{ZhMe2016} the closures of the invariant manifolds of the saddles are locally flat spheres. That is a contradiction because the closures of the invariant manifolds of the saddle singularities of $X^t_0$ and $f_0$ coincide. Thus, $f_0$ does not embed into a flow.

\begin{figure}[h]
\center{\includegraphics[width=0.6\linewidth]{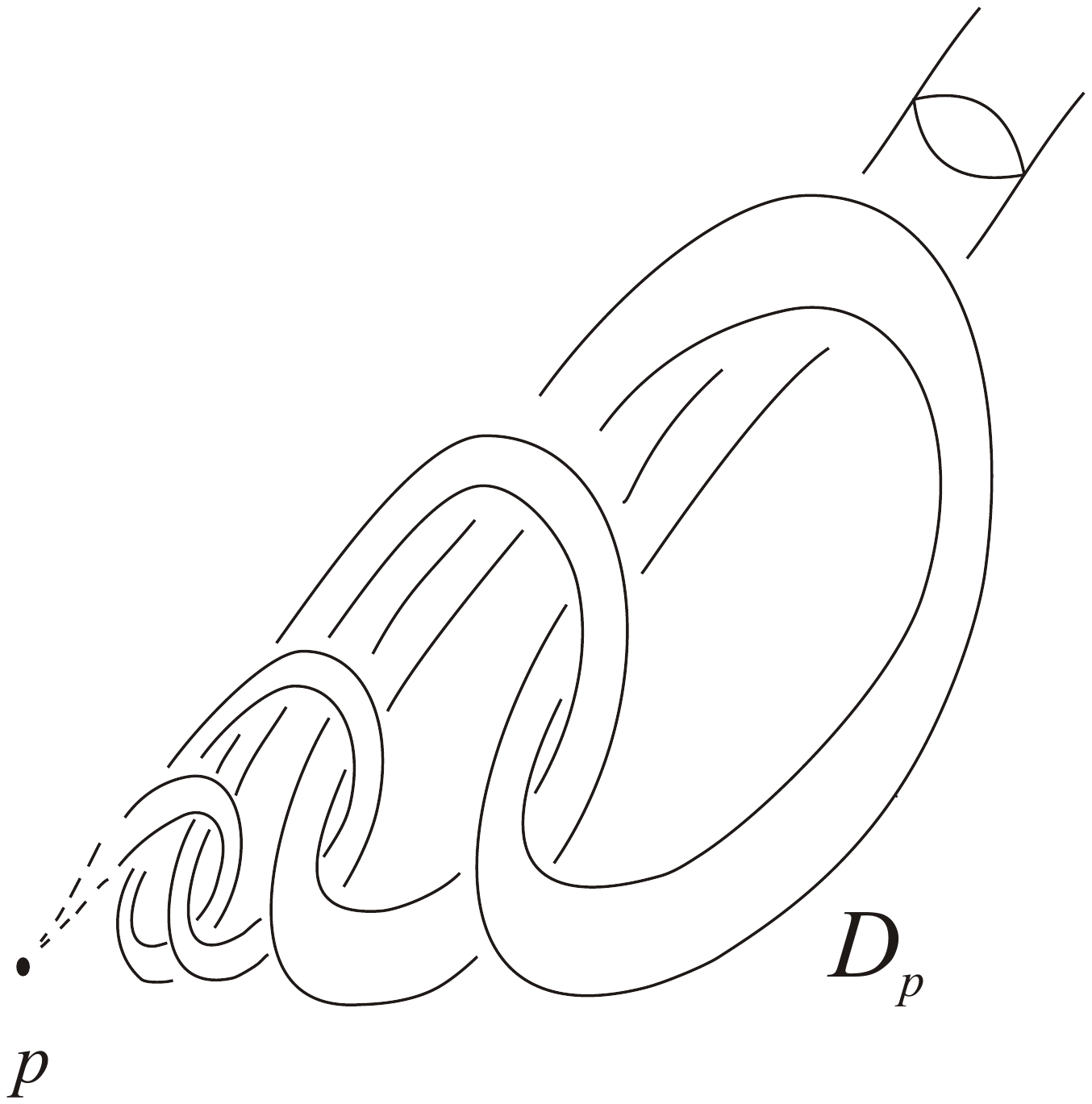}}
\caption{The disk $D_p\subset W^s_p$}
\label{wild-disk}
\end{figure}

 In~\cite{MePo} T.Medvedev and O. Pochinka constructed  an example of Morse-Smale diffeomorphism  $f_1:S^4\to S^4$ satisfying to the conditions  $i)-ii)$ of the Theorem~\ref{emb-to-flow}.  The non-wandering set of the diffeomorphism  $f_1$ consists of two sources, two sinks and two saddles  $p,q$  such that $dim~W^s_p=dim~W^u_q=3$. The  intersection  $W^s_p\cap W^u_q$  is not empty and its closure in  $W^s_p$ is a wildly embedded open disk $D_p$ (see Fig.~\ref{wild-disk}). If $S^2\subset W^s_p$ is a  2-sphere   which bounds an open  ball containing  the point  $p$  then the intersection  $S^2\cap D_p$ contains at least three connected components. Assume that $f_1$ embeds into a topological flow $X^t_1$. Then due to \cite{GrGuMePo}  the restriction $X^t_p$ of $X^t_1$ to $W^s_p\setminus p$  is topologically conjugated by means of a homeomorphism $h:W^s_p\setminus p\to\mathbb S^2\times\mathbb R$ to a shift flow $\chi^t(s,r)=(s,r+t),\,(s,r)\in\mathbb S^2\times\mathbb R$. Let $\Sigma=h^{-1}(\mathbb S^2\times\{0\})$. Then every trajectory of the flow  $X^t_p$  intersects the sphere $\Sigma$  at  a unique point. Since the disk  $D_p$ is invariant with respect to the flow  $X^t_p$ the intersection  $D_p\cap \Sigma$ consists of a unique connected component and that is a contradiction.  Thus, $f_1$ does not embed into a flow.

\section{The scheme of the proof of  Theorem~\ref{emb-to-flow}}\label{sort}

The proof  of  Theorem~\ref{emb-to-flow} is based on the technique developed for classification of Morse-Smale diffeomorphisms on orientable manifolds in a series of papers~\cite{BoGr00}, \cite{BoGrMePe}, \cite{BoGrPo05}, \cite{GrPo}, \cite{GrGu07}, \cite{GrGuMe08}, \cite{GrGuMe10},\cite{GrGuPo-matan}.  The idea of the proof consists of the following. 

In section~\ref{MShom} we introduce a  notion of Morse-Smale homeomorphism on a topological n-manifold and define  the subclass $G(S^n)$  of such  homeomorphisms satisfying to conditions similar to $i)-iii)$  of  Theorem~\ref{emb-to-flow}. 

Let $f\in G(S^n)$. In~\cite[Theorem 1.3]{GrGuPo-matan} it is shown that the dimension of the invariant manifolds of the fixed points of $f$ can be only one of  $0,1,n-1$ or $n$. Denote by $\Omega^{i}_f$ the set of all fixed points of $f$ whose unstable manifolds have dimension  $i\in \{0,1,n-1,n\}$, and  by $m_f$ the number of all saddle points of $f$. 

Represent the sphere   $S^n$  as the union of pairwise disjoint sets  $$A_f=(\bigcup\limits_{\sigma\in \Omega^{1}_f}{W^u_\sigma})\cup \Omega^{0}_f,\, R_f=(\bigcup\limits_{\sigma\in \Omega^{n-1}_f}{W^s_\sigma})\cup  \Omega^{n}_f,\,V_f=S^n\setminus(A_f\cup R_f).$$  Similar to~\cite{GrPoZh} one can prove that the sets  $A_f, R_f, V_f$  are connected, the set   $A_f$  is an attractor, $R_f$ is a repeller\footnote{A set  $A$ is called an attractor of a homeomorphism  $f:M^n\to M^n$ if there exists a closed neighborhood  $U\subset M^n$ of the set  $A$ such that $f(U)\subset int~U$ and  $A=\bigcap\limits_{n\geq 0} f(U)$.  A set  $R$ is called a repeller of a homeomorphism  $f$ if it is an attractor for the homeomorphism $f^{-1}$.} and $V_f$ consists of wandering orbits of $f$ moving from  $R_f$ to $A_f$. 

Denote by  $\widehat V_f=V_f/f$  the orbit space of the action of  $f$ on $V_f$ and by   $p_{_f}:V_f\to \widehat V_f$ the natural projection. Let $$\hat{L}^s_{f}=\bigcup \limits_{\sigma\in \Omega^1_f} p_{_f}(W^s_\sigma\setminus\sigma),\,\,\, \hat{L}^u_{f}=\bigcup \limits_{\sigma\in \Omega^{n-1}_f} p_{_f}(W^u_\sigma\setminus\sigma).$$ 

\begin{Defin} The collection  $S_{f}=(\widehat V_{f},\hat{L}^s_{f},\hat{L}^u_{f})$ is called the scheme of the homeomorphism  $f\in G(S^n)$.
\end{Defin}

\begin{Defin}  Schemes $S_f$ and  $S_{f'}$ of homeomorphisms $f,f'\in G(S^n)$ are called equivalent if there exists a homeomorphism  $\hat\varphi:\widehat V_f\to\widehat V_{f'}$ such that   $\hat\varphi(\hat{L}^s_{f})=\hat{L}^s_{f'}$ and  $\hat\varphi(\hat{L}^u_{f})=\hat{L}^u_{f'}$. 
\end{Defin}

The next statement follows from  paper~\cite[Theorem 1.2]{GrGuPo-matan} (in fact, Theorem 1.2 was proven for Morse-Smale diffeomorphisms but the smoothness plays no role in the proof).

\begin{st}\label{mat} Homeomorphisms  $f,f'\in G(S^n)$ are topologically equivalent if and only if their schemes $S_f$,  $S_{f'}$ are equivalent.
\end{st}

The possibility of embedding of $f\in G(S^n)$ into a topological flow follows from triviality of the scheme in the following sense.

Let $a^t$ be the flow on the set  $\mathbb S^{n-1}\times \mathbb{R}$ defined by $a^t(x,s)=(x,s+t)$, $x\in S^{n-1}, s\in \mathbb{R}$ and let $a$ be the time-one map of $a^t$. Let $\mathbb{Q}^n=\mathbb{S}^{n-1}\times \mathbb{S}^1$. Then the orbit space of the action  $a$ on $\mathbb{S}^{n-1}\times \mathbb{R}$ is $\mathbb Q^n$. Denote by  $p_{_{\mathbb Q^n}}:\mathbb{S}^{n-1}\times \mathbb{R}\to\mathbb Q^n$  the natural projection.
Let $m\in\mathbb N$ and  $c_1,...,c_m\subset\mathbb S^{n-1}$ be a collection of smooth pairwise disjoint $(n-2)$-spheres. Let ${Q}^{n-1}_i=\bigcup \limits_{t\in \mathbb{R}} a^t(c_i)$, $\mathbb L_m=\bigcup \limits_{i=1}^{m}Q^{n-1}_i$ and $\widehat{\mathbb L}_m=p_{_{\mathbb Q^n}}(\mathbb L_m)$.

\begin{Defin}\label{simple}  The scheme  $S_{f}=(\widehat V_{f},\hat{L}^s_{f},\hat{L}^u_{f})$  of a homeomorphism  $f\in G(S^n)$ is called trivial if there exists a homeomorphism $\hat\psi:\widehat V_f\to\mathbb Q^n$  such that   $\hat\psi(\hat{L}^s_{f}\cup \hat{L}^u_f)=\widehat{\mathbb L}_{m_f}$. 
\end{Defin}

In the section~\ref{triv-pr}  we prove the following key lemma.

\begin{Lemm}\label{simple}  If $f\in G(S^n)$ then its scheme $S_{f}$ is trivial.
\end{Lemm}   

In the section~\ref{flow} we construct a topological flow $X^t_f$  whose  time one map belongs to the class  $G(S^n)$ and has the  scheme  equivalent to $S_f$. According to Statement~\ref{mat}  there exists a homeomorphism  $h:S^n\to S^n$ such that  $f=hX^1_fh^{-1}$. Then the homeomorphism  $f$ embeds into the topological flow  $Y^t_f=hX^t_fh^{-1}$.  
\section{{Morse-Smale homeomorphisms}}
\label{MShom}

This section contains   some definitions and statements which was introduced and proved in \cite{GrGuMePo-mzm2017}.
\subsection{Basic definitions}

Remind that a linear automorphism  $L:\mathbb{R}^n\to \mathbb{R}^n$  is called hyperbolic if its matrix has no eigenvalues with absolute value equal one. In this case a space  $\mathbb{R}^n$  have a unique  decomposition into the direct sum of   $L$-invariant subsets  $E^s, E^u$  such that  $||L|_{E^s}||<1$ and $||L^{-1}|_{E^u}||< 1$ in some norm  $||\cdot||$ (see, for example, Propositions 2.9, 2.10 of Chapter 2 in~\cite{PaMe}).   

According to Proposition~5.4 of the book~\cite{PaMe} any hyperbolic automorphism   $L$ is topologically conjugated with a  linear map of the following form:

\begin{equation}\label{map}
a_{\lambda,\mu,\nu}(x_1,x_2,...,x_{\lambda}, x_{\lambda+1},x_{\lambda+2},...,x_n)=(2\mu x_1,2x_2,...,2x_{\lambda},\frac12 \nu x_{\lambda+1},\frac12 x_{\lambda+2},...,\frac12 x_n),
\end{equation}

 where  $\lambda=dim\, E^u\in\{0,1,...,n\}$,   $\mu=-1$ ($\mu=1$)  if the restriction  $L|_{E^u}$  reverses (preserves)  an orientation of  $E^u$,  and  $\nu=-1$ ($\nu=1$) if the restriction  $L|_{E^s}$ reverses (preserves) an orientation of $E^s$.  

Put ${\mathbb{E}^s_\lambda}=\{(x_1,...,x_n)\in \mathbb{R}^n|\ x_1=x_2=\dots=x_{\lambda}=0\}$, ${\mathbb{E}^u_\lambda}=\{(x_1,...,x_n)\in \mathbb{R}^n|\ x_{\lambda+1}=x_{\lambda+2}=\dots=x_{n}=0\}$ and denote by  $P^s_x (P^u_y)$ a hyperplane that  parallel to the hyperplane  $\mathbb{E}^s_\lambda$ ($\mathbb{E}^u_\lambda$) and contain a point  $x\in \mathbb{E}^u_\lambda$ ($y\in \mathbb{E}^s_\lambda$). Unions $\mathcal{P}^s_\lambda=\{P^s_x\}_{x\in \mathbb{E}^u_\lambda}, \mathcal{P}^u_\lambda=\{P^u_y\}_{y\in \mathbb{E}^s_\lambda}$  form  the  $a_{\lambda,\mu,\nu}$-invariant foliation.

Suppose that  $M^n$ is an $n$-dimensional topological manifold,  $f:M^n\to M^n$ is a homeomorphism  and   $p$ is a fixed point of the homeomorphism  $f$. We will call the point $p$ {\it topologically hyperbolic point of index  $\lambda_p$}, if there exists its neighborhood  $U_p\subset M^n$, numbers  $\lambda_p\in\{0,1,...,n\},\mu_p,\nu_p\in\{+1,-1\}$, and a homeomorphism  $h_p:U_p\to \mathbb{R}^n$ such that $h_pf|_{U_p}=a_{\lambda_p,\mu_p,\nu_p}h_p|_{U_p}$   when the left and right parts are defined.  Call the  sets $W^s_{p,loc}=h_p^{-1}(E^s), {W}^u_{p,loc}= h_p^{-1}(E^u)$     {\it the  local invariant manifolds} of the point $p$, and the   sets  $W^s_p=\bigcup \limits_{i\in \mathbb{Z}}f^i({W}^s_{p,loc})$, $W^u_p=\bigcup \limits_{i\in \mathbb{Z}}f^i({W}^u_{p,loc})$     {\it the stable and unstable invariant manifolds of the point $p$}.  

It follows form the definition that   $W^s_p=\{x\in M^n: \lim\limits_{i\to +\infty} f^{i}(x)= p\}, W^u_p=\{x\in M^n: \lim\limits_{i\to +\infty} f^{-i}(x)= p\}$ and  $W^u_p\cap W^u_q=\emptyset$ ($W^s_p\cap W^s_q=\emptyset$) for any distinct hyperbolic points $p,q$. Moreover,  there exists  an injective  continuous immersion  $J:\mathbb{R}^{\lambda_p}\to M^n$ such that  $W^u_p=J(\mathbb{R}^{{\lambda_p}})$\footnote{A map $J:\mathbb{R}^{m}\to M^n$ is called immersion if for any point  $x\in \mathbb{R}^{m}$  there exists a neighborhood  $U_x\in \mathbb{R}^{m}$  such that the restriction  $J|_{U_x}$ of the map $J$ on the set  $U_x$ is a homeomorphism.}.  

A hyperbolic fixed   point is called  {\it the source} ({\it the sinks})  if its  indice equals   $n$ ($0$), a hyperbolic fixed   point $p$ of index  $0<\lambda_p<n$  is called   {\it the saddle point}.

A periodic point  $p$ of  period  $m_p$ of a homeomorphism  $f$ is called     {\it  a topologically hyperbolic  sink $($ source, saddle$)$ periodic point} if it is   the  topologically hyperbolic $($source, saddle$)$ fixed point for the homeomorphism  $f^{m_p}$. The  stable and unstable manifolds of the periodic  point  $p$   considered  as the fixed point of the homeomorphism  $f^{m_p}$  are called the    stable and unstable manifolds of the point  $p$.  Every connected component of the set   $W^s_{p}\setminus p$ ($W^u_{p}\setminus p$)  is called  {\it the  stable  $($ the  unstable$)$ separatrix}  and is denoted  by  $l^s_{p}$ ($l^u_{p}$). 

 The linearizing   homeomorphism  $h_p:U_p\to \mathbb{R}^n$ induces   a pair of transversal foliations  $\mathcal{F}^s_p=h_p^{-1}(\mathcal{P}^s_{\lambda_p})$ , $\mathcal{F}^u_p=h_p^{-1}(\mathcal{P}^u_{\lambda_p})$ on the set  $U_p$.  Every leaf of the foliation $\mathcal{F}^s_p\  (\mathcal{F}^u_p)$  is an open  disk of dimension $\lambda_p\   (n-\lambda_p)$. For any point $x\in U_p$ denote by  $F^s_{p,x},F^u_{p,x}$  the  leaf of the foliation  $\mathcal{F}^s_p,\mathcal{F}^u_p$, correspondingly,  containing the point $x$.

The invariant manifolds   $W^s_p$ and $W^u_q$ of saddle periodic points $p,q$ of a homeomorphism  $f$  intersect  {\it consistently transversally} if one of the following conditions holds:
\begin{enumerate}
\item $W^s_p\cap W^u_q=\emptyset$;
\item $W^s_p\cap W^u_q\neq\emptyset$ and $F^s_{q,x}\subset W^s_p$; $F^u_{p,y}\subset W^u_q$ for any points $x\in W^s_p\cap U_q$,  $y\in W^u_q\cap U_p$.
\end{enumerate}

\begin{Defin}  A homeomorphism $f:M^n\to M^n$  is called the Morse-Smale homeomorphism if it satisfies  the next conditions:
\begin{enumerate}
\item its  non-wandering set  $\Omega_f$ finite and any point  $p\in\Omega_f$  is topologically hyperbolic;
\item invariant manifolds of any two saddle points $p,q\in \Omega_f$  intersect   consistently  transversally.
\end{enumerate}
\end{Defin}

\subsection{Properties of Morse-Smale homeomorphisms}\label{MS-properties}

\begin{st}\label{homo}  Let  $f:M^n\to M^n$ be a Morse-Smale homeomorphism. Then:
\begin{enumerate}
\item $W^u_p\cap W^s_p=p$ for any saddle point  $p\in \Omega_f$;
\item for any saddle points  $p,q,r\in \Omega_f$ the conditions $(W^s_p\setminus p)\cap (W^u_q\setminus q)\neq \emptyset$, $(W^s_q\setminus q)\cap (W^u_r\setminus r)\neq \emptyset$ imply  $(W^s_p\setminus p)\cap (W^u_r\setminus r)\neq \emptyset$;
\item  there are no sequence of distinct saddle points $p_1,p_2,...,p_{k}\in \Omega_f$, $k>1$, such that  $(W^s_{p_i}\setminus p_i)\cap (W^u_{p_{i+1}}\setminus p_{i+1})\neq \emptyset$ for $i\in \{1,...,k-1\}$ and $(W^s_{p_k}\setminus p_k)\cap (W^u_{p_1}\setminus p_1)\neq \emptyset$.
\end{enumerate}
\end{st}

\begin{st}\label{smale} Let  $f:M^n\to M^n$ be a Morse-Smale homeomorphism. Then: 
\begin{enumerate}
\item[$1)$] $M^n=\bigcup\limits_{p\in \Omega_f}W^u_p$;
\item[$2)$] for any point  $p\in \Omega_f$  the manifold  $W^u_p$  is a topological submanifold of the manifold  $M^n$;
\item[$3)$] for any point    $p\in \Omega_f$ and any connected component  $l^u_p$  of the set $W^u_p\setminus p$ the following equality holds:  $cl\,{l^{u}_p}\setminus(l^{u}_p\cup
p)=\bigcup\limits_{q\in \Omega_f:W^s_q\cap l^{u}_p\neq \emptyset}W^u_q$ \footnote{Here  $cl\,{l^{u}_p}$ means the closure of the set  $l^{u}_p$.}.
\end{enumerate}
\end{st}

\begin{Cor}\label{1-sink-cor}
If   $f:M^n\to M^n$ is a Morse-Smale homeomorphism and  $p\in \Omega_f$ is a  saddle point such that  $l^u_{p}\cap W^s_q=\emptyset$   for any saddle point  $q\neq p$, then there exists a unique sink  $\omega\in \Omega_f$ such that $cl\,l^{u}_p=l^u_p\cup p\cup \omega$ and   $cl\,l^{u}_p$ is either  a compact arc in case  $\lambda_p=1$ or a  sphere of dimension $\lambda_p$  in case  $\lambda_p>1$.
\end{Cor}

For an arbitrary point $q\in \Omega_f$ and $\delta\in\{u,s\}$ put $V^\delta_q=W^\delta_q\setminus q$ and denote by  $\widehat{V}^\delta_q=V^s_q/f$ the orbit space of the action of the homeomorphism $f$  on the set $V^\delta_q$. The following statement is proved in the book~\cite{GrPo} (Proposition~2.1.5).

\begin{st}\label{smoothconj-cor1} The space  $\widehat{V}^u_q$ is homeomorphic  to  $\mathbb{S}^{\lambda_q-1}\times \mathbb{S}^{1}$ and the space  $\widehat{V}^s_q$ is homeomorphic to  $\mathbb{S}^{n-\lambda_q-1}\times \mathbb{S}^{1}$.
\end{st}

Remark  that   $\mathbb{S}^{0}\times \mathbb{S}^{1}$ means   a union of two disjoint  closed curves.

\begin{Prop} \label{cil} 
Suppose  $f:M^n\to M^n$  is a Morse-Smale homeomorphism, $n\geq 4$, and  $\sigma \in \Omega_f$ is a saddle point  of   index $(n-1)$ such that  $l^u_{\sigma}\cap W^s_q=\emptyset$ for any saddle point $q\neq p$.  Then the sphere  $cl\,{l^u_\sigma}$ is bicollared.
\end{Prop}
\begin{demo}
Let    $\omega\in\Omega^0_f$ be a sink point such that $l^u_\sigma\subset W^s_\omega$.  
Due to  Corollary~\ref{1-sink-cor} and the item  2 of  Statement~\ref{smale} the set    $cl\,{l^u_\sigma}=l^u_\sigma\cup \omega$ is an  $(n-1)$-sphere which is locally flat embedded in $M^n$ at all its points apart possibly one point $\omega$. According to~\cite{Ca63}, \cite{Ki68}  an $(n-1)$-sphere in a manifold  $M^n$ of dimension  $n\geq 4$ is  either locally flat or have more than countable set of  points of wildness.  Therefore the sphere  $cl\,{l^u_\sigma}$ is locally flat at point  $\omega$. According to~\cite{Br62} a locally flat sphere is bicollared.
\end{demo}

By  $G(S^n)$  we denoted a class of Morse-Smale homeomorphism on the sphere $S^n$ such that  any  $f\in G(S^n)$  satisfy the following conditions:
\begin{enumerate}
\item[$i)$] $\Omega_f$ consists of fixed points;
\item[$ii)$] $W^s_p\cap W^u_q=\emptyset$ for any distinct saddle points  $p,q\in \Omega_f$;
\item[$iii)$] the restriction of a homeomorphism $f$ on every invariant manifolds of an arbitrary fixed point  $p\in \Omega_f$ preserves its orientation.
\end{enumerate}

\begin{Prop}\label{sphere}
If $f\in G(S^n)$, then any  saddle fixed point has index   $1$ and $(n-1)$.
\end{Prop}
\begin{demo}
Suppose that, on the contrary, there exists a point  $\sigma\in \Omega_f$ of index   $j\in (1, n-1)$. According to  Corollary~\ref{1-sink-cor} the closures 
 $cl\,{W^{u}_{\sigma}}, cl\,{W^{s}_{\sigma}}$ of the stable  and unstable manifolds of the point  $\sigma$ are spheres of dimensions $j$ and $n-j$ correspondingly. Due to   item 1 of  Statements~\ref{cil},  the spheres $S^j=cl\,{W^{u}_\sigma},S^{n-j}= cl\,{W^{s}_\sigma}$   intersect at a single point  $\sigma$. Therefore their intersection index equals either  $1$  or $-1$ (depending on the choice of orientations of the spheres $S^j$, $S^{n-j}$ and $S^n$). Since homology groups  $H_j(S^n), H_{n-j}(S^n)$ are trivial it follows that  there is a sphere  $\tilde{S}^{j}$  homological to the sphere  $S^j$ and having the empty intersection with the sphere   $S^{n-j}$.  Then the intersection number of the spheres    $S^j, S^{n-j}$ must be equal to zero as the  intersection number is the homology invariant  (see, for example, \cite{ZeTr}, $\S\,69$). This contradiction proves the statement.    
\end{demo}

\subsection{Canonical manifolds connected with saddle fixed points of a homeomorphism  $f\in G(S^n)$}\label{aux}
\hspace*{\parindent} 

It follows from  Statement~\ref{sphere}  that for each  saddle point  of a  homeomorphism  $f\in G(S^n)$ there exists a neighborhood where $f$  is topologically conjugated either with the map  $a_1:\mathbb{R}^{n}\to \mathbb{R}^n$ defined by  $a_1(x_1,x_2,\dots,x_n)=(2 x_1,\frac12 x_{2},\dots,\frac12 x_n)$ or with the map $a_1^{-1}$.  In this section we describe canonical manifolds defined by  the action of the map $a_1$ and prove  Proposition~\ref{grgupo-ado} allowing to define similar canonical manifolds for the homeomorphism $f\in G(S^n)$.

Put $\mathbb{U}_{\tau}=\{(x_{1},...,x_{n})\in
\mathbb{R}^{n}|\ x_{1}^{2}(x_{2}^{2}+...+x_{n}^{2})\leq \tau^2\}$, $\tau\in (0,1]$, $\mathbb{U}=\mathbb{U}_1$; $\mathbb{U}_0=\{(x_1,...,x_n)\in \mathbb{R}^n|\ x_1=0\}$, $\mathbb{N}^s=\mathbb{U}\setminus Ox_1$,  $\mathbb{N}^u=\mathbb{U}\setminus \mathbb{U}_0$,  $\widehat{\mathbb{N}}^{s}=\mathbb{N}^s/a_1$, $\widehat{\mathbb{N}}^{u}=\mathbb{N}^u/a_1$.  Denote by  $p_s: \mathbb{N}^s\to \widehat{\mathbb{N}}^{s}$, $p_u: \mathbb{N}^u\to \widehat{\mathbb{N}}^{u}$ the  natural projections and put  $\widehat{V}^{s}=p_{s}(\mathbb{U}_0)$.

 \begin{figure}[h]
\center{\includegraphics[width=0.8\linewidth]{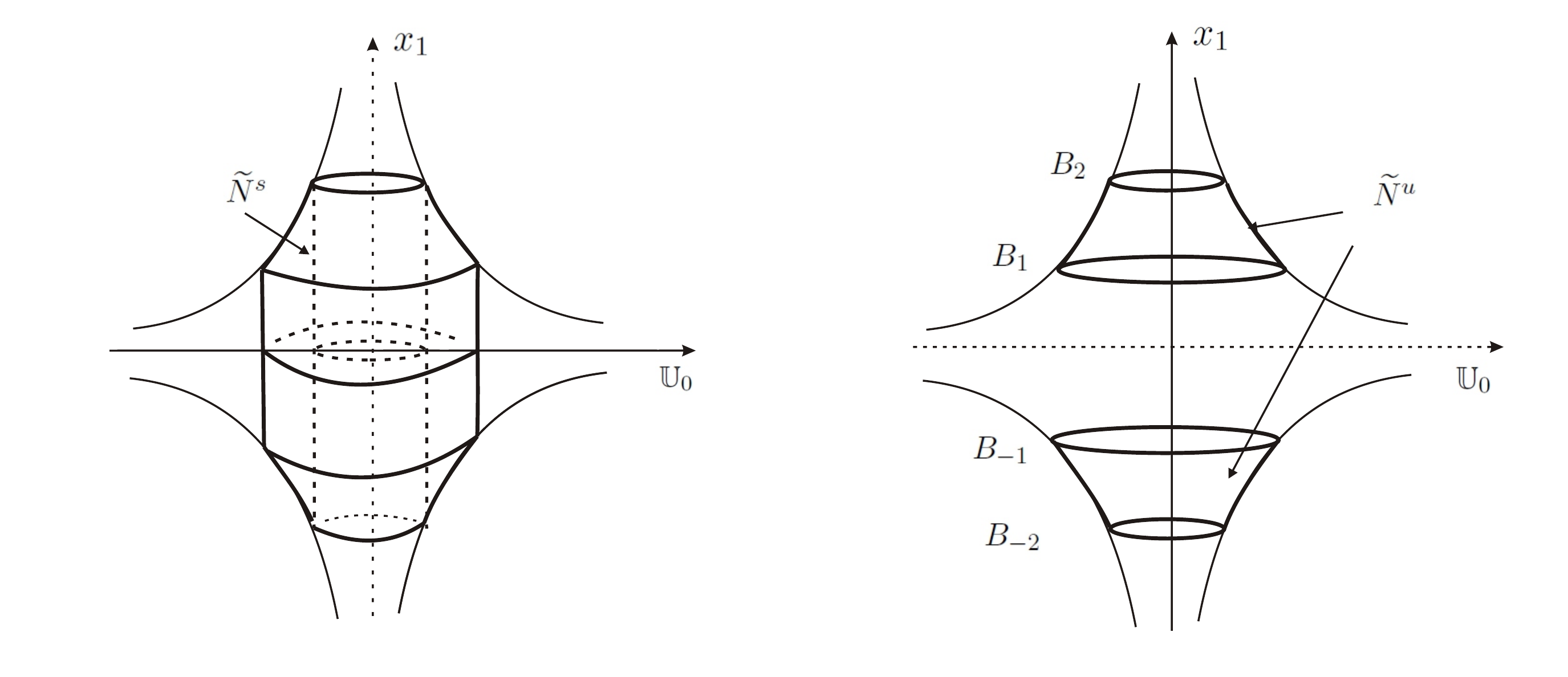}}
\caption{Fundamental domains  $\widetilde{N}^s, \widetilde{N}^u$ of the action of the homeomorphism  $a_1$ on the sets  $\mathbb{N}^s, \mathbb{N}^u$}
\label{canon-nbh}
\end{figure}

 The following statement is proved in~\cite{GrGuMe10} (Propositions 2.2, 2.3).

\begin{st}\label{st-nbh} 
The space  $\widehat{\mathbb{N}}^{s}$ is homeomorphic to the direct product $\mathbb{S}^{n-2}\times \mathbb{S}^1\times [-1,1]$,  the space  $\widehat{\mathbb{N}}^{u}$  consists of two connected components each of which is homeomorphic to the direct product   $\mathbb{B}^{n-1}\times \mathbb{S}^1$.
\end{st}

Recall that  {\it an annulus} of dimension  $n$  is a manifold homeomorphic to $\mathbb{S}^{n-1}\times [0,1]$. 

{ On the Figure~\ref{canon-nbh} we present the neighborhoods  $\mathbb{N}^s, \mathbb{N}^u$ and the fundamental domains $\widetilde{N}^s=\{(x_1,\dots,x_n)\in \mathbb{N}^s|  \frac14\leq x_2^2+\dots+x_n^2\leq 1\}, \widetilde{N}^u=\{(x_1,\dots,x_n)\in \mathbb{N}^u| |x_1|\in [1,2]\}$ of the action of the diffeomorphism  $a_1$\footnote{A fundamental domain of the action of a group  $G$ on a set  $X$ is a closed set  $D_{_{G}}\subset X$ containing a subset  $\tilde{D}_{_{G}}$  with the following properties: 1) $cl\,\tilde{D}_{_{G}}=D_{_{G}}$; 2) $g(\tilde{D}_{_{G}})\cap \tilde{D}_{_{G}}=\emptyset$  for any  $g\in G$ distinct  from the neutral element; 3) $\bigcup\limits_{g\in G}g(\tilde{D}_{_{G}})=X$.}. Put   $\mathcal{C}=\{\{(x_1,\dots,x_n)\in \mathbb{R}^n|  \frac14\leq x_2^2+\dots+x_n^2\leq 1\}$.  The set  $\mathbb{N}^s$  is the union of  the hyperplanes  $\mathcal{L}_t=\{(x_1,...,x_n)\in {N}^s| x_1^2(x_2^2+\dots+x_n^2)=t^2\}, t\in [-1,1]$. Then the fundamental domain  $\widetilde{{N}}^{s}$ is the   union of the pairs of annuli   $\mathcal{K}_t=\mathcal{L}_t\cap \mathcal{C}, t\in [-1,1]$  and the space   $\widehat{\mathbb{N}}^{s}$ can be obtained from  $\widetilde{{N}}^{s}$  by gluing the  connected components of the boundary of each annulus by means of the  diffeomorphism  $a_1$. The set  $\widetilde{N}^u$ consist of two connected components each of which is homeomorphic to the direct product  $\mathbb{B}^{n-1}\times [0,1]$.  The space   $\widehat{\mathbb{N}}^{u}$   is obtained  from $\widetilde{N}^u$ by gluing the disk   $B_1=\{(x_1,\dots,x_n)\in \mathbb{N}^u| x_1=1\}$  to the disk  $B_2=\{(x_1,\dots,x_n)\in \mathbb{N}^u| x_1=2\}$ and the disk  $B_{-1}=\{(x_1,\dots,x_n)\in \mathbb{N}^u| x_1=-1\}$  to the disk $B_{-2}=\{(x_1,\dots,x_n)\in \mathbb{N}^u| x_1=-2\}$ by means of the diffeomorphism  $a_1$.}

\begin{Prop}\label{grgupo-ado} Suppose $f\in G(S^n)$; then there exists a set of pair-vise disjoint neighborhoods  $\{N_\sigma\}_{\sigma\in \Omega^1_f\cup \Omega^{n-1}_f}$   such that for any neighborhood $N_{\sigma}$  there exists a homeomorphism  $\chi_{\sigma}:N_{\sigma} \to \mathbb{U}$ such that  $\chi_{\sigma}f|_{_{N_{\sigma}}}=a_{1}\chi_{\sigma}|_{_{N_{\sigma}}}$ whenever $\lambda_\sigma=1$ and  $\chi_{\sigma}f|_{_{N_{\sigma}}}=a^{-1}_{1}\chi_{\sigma}|_{_{N_{\sigma}}}$ whenever $\lambda_\sigma=n-1$.
\end{Prop}
\begin{demo}
Put   ${V}^\delta_{\Omega^i_f}=\bigcup\limits_{q\in \Omega^i_f}{V}^\delta_q$, $\widehat{V}^\delta_{\Omega^i_f}=\bigcup\limits_{q\in \Omega^i_f}\widehat{V}^\delta_q$, $i\in \{0,1, n-1, n\}$, $\delta\in \{s,u\}$ and denote by $p^{\delta}_{\Omega^i_f}:{V}^\delta_{\Omega^i_f}\to \widehat{V}^\delta_{\Omega^i_f}$ the natural projection such that  $p^{\delta}_{\Omega^i_f}|_{V^\delta_q}=p^{\delta}_q|_{V^\delta_q}$ for any point   $q\in \Omega_f$.
	
	Put $\Sigma_f=\Omega^1_f\cup \Omega^{n-1}_f$,  $\widehat{L}^u_{\Sigma_f}=p^{s}_{\Omega^0_f}(V^u_{\Omega^1_f}\cup V^u_{\Omega^{n-1}_f})$.

	
	The set  $\widehat{L}^u_{\Sigma_f}$ consists of  finite number of compact topological submanifolds. Then  
	there is a  set of pair-vise disjoint compact neighborhoods   $\{\widehat{K}^u_\sigma,\sigma\in  \Sigma_f\}$ 
	of these manifolds in  $\widehat{V}^s_{\Omega^0}$. For every point  $\sigma\in \Sigma_f$ put $K^u_\sigma=(p^s_{\Omega^0_f})^{-1}({\widehat{K}^u_\sigma})$ and   $\widetilde{N}_\sigma=K^u_\sigma\cup W^s_\sigma$. 

 Let   $U_\sigma\subset \widetilde{N}_\sigma$ be a neighborhood of the point   $\sigma$ such that  a homeomorphism   $g_\sigma:U_\sigma\to \mathbb{R}^n$ satisfying the condition   $g_\sigma f|_{U_\sigma}=a_{\lambda_\sigma} g_\sigma|_{U_\sigma}$ is defined. 

   Put $u_{\tau}=\{(x_1,...,x_n)\in \mathbb{U}_\tau|\ x_2^2+...+x_{n}^2\leq 1, |x_1|\leq 2\tau\}$, $D^u_\tau=\{(x_1,...,x_n)\in \mathbb{U}_\tau|\  \tau<|x_1|\leq 2\tau \}$, $D^s_\tau=\{(x_1,...,x_n)\in \mathbb{U}_\tau|\  \frac14\leq x_2^2+...+x_{n}^2\leq 1\}$, $\tilde{u}_\tau=g_{\sigma}^{-1}(u_\tau)$, $\widetilde{D}^\delta_\tau=g_{\sigma}^{-1}(D^\delta_\tau)$, $\delta\in \{s,u\}$, and  $N_\tau=\bigcup \limits_{i\in \mathbb{Z}}f^i(\widetilde{u}_\tau)$.
	
	Let us show that there is a  number $\tau_1>0$  such that for any $i\in \mathbb{N}$ the intersection $f^i(\widetilde{D}^u_{\tau_1})\cap \tilde{u}_{\tau_1}$ is empty.   Suppose  $\sigma\in \Omega^{n-1}_f$ (the argument  for the case $\sigma\in \Omega^1_f$ is similar). By the Statement~\ref{smale},  the set $\bigcup\limits_{i\in \mathbb{N}}f^i (\widetilde{D}^u_\tau)$ lies in the stable manifold of a unique sink  point $\omega$.  Since the homeomorphism $f$ is locally conjugated with the  linear compression $a_{0}$ in a neighborhood of the point $\omega$, we have  that  there exists  a ball $B^n\subset W^s_\omega\setminus U_\sigma$ such that $\omega\subset B^n$ and $f(B^n)\subset int~B^{n}$. Since  $\widetilde{D}^u_\tau$ is compact, there is $i^*>0$ such that $f^{i}(\widetilde{D}^u_\tau)\cap U_\sigma\subset B^n$ for all $i>i^*$. Hence the set of numbers $i_j$ such that $f^{i_j}(\widetilde{D}^u_\tau)\cap \tilde{u}_\tau\neq \emptyset$ is finite.   Then one can choose $\tau_1\in (0,\tau)$ such that $\tilde{u}_{\tau_1}\cap f^i(\widetilde{D}^u_\tau))=\emptyset$ and therefore $\tilde{u}_{\tau_1}\cap f^i(\widetilde{D}^u_{\tau_1}))=\emptyset$ for any $i\in \mathbb{N}$.  Similarly one can show  that there exists   a number  $\tau_2\in (0,\tau_1]$ such that for any $i\in \mathbb{N}$ the intersection of $f^{-i}(\widetilde{D}^s_{\tau_2})\cap \tilde{u}_{\tau_2}$ is empty.
				
			Suppose  $\lambda_\sigma=1$, put $N_\sigma=\bigcup\limits_{i\in \mathbb{Z}}f^i(\tilde{u}_{\tau_2})$,  and define a homeomorphism  ${\chi^*}_\sigma: N_\sigma\to U_{\tau_2}$  by  the following:  $\chi^*_\sigma(x)=g_\sigma(x)$ whenever  $x\in \tilde{u}_{\tau_2}$,   and   $\chi^*_\sigma(x)=a_{\lambda_\sigma}^{-k}(g_\sigma(f^{k}(x)))$ whenever  $x\in {N}_\sigma\setminus (\tilde{u}_{\tau_2})$, where $k\in \mathbb{Z}$ is  such that  $f^k(x)\in \tilde{u}_{\tau_2}$. The homeomorphism  ${\chi^*}_\sigma$  conjugates the   homeomorphism  $f|_{N_\sigma}$ with the linear diffeomorphism  $a_{1}|_{\mathbb{U}_{\tau_2}}$.  Since the homeomorphism   $a_{1}|_{\mathbb{U}_{\tau_2}}$ is topologically conjugated with $a_{1}|_{\mathbb{U}}$ by means of  the diffeomorphism  $g(x_1,...,x_n)=\left(\frac{x_1}{\sqrt{\tau_2}},...,\frac{x_n}{\sqrt{\tau_2}}\right)$, we see that  the superposition   $\chi_\sigma=g\chi^*_\sigma:{N}_\sigma\to \mathbb{U}$ topologically conjugates   $f|_{{N}_\sigma}$ with $a_{1}|_{\mathbb{U}}$.  A homeomorphism   $\chi_\sigma$  for the case  $\lambda_\sigma=n-1$ can be constructed in the same way.

\end{demo}

Put $N^u_\sigma=N_\sigma\setminus W^s_\sigma$, $N_{\tau,\sigma}=\chi_\sigma^{-1}(\mathbb{U}_\tau)$, $N^s_\sigma=N_\sigma\setminus W^u_\sigma$, $\widehat{N}^{s}_\sigma=N^s_\sigma/f$, $\widehat{N}^{u}_\sigma=N^u/f$.

\section{Triviality of the scheme of the homeomorphism  $f\in G(S^n)$}\label{triv-pr}

This section is devoted to the proof of   Lemma~\ref{simple}. In subsections~\ref{1-emb}-\ref{inv-surg} we establish some  axillary results.  

\subsection{Introduction results on the embedding of closed curves and their tubular neighborhoods in a manifold $M^n$}\label{1-emb}

Further we  denote by  $M^n$ a topological manifold possibly with non-empty boundary. 
 
 Recall that  a manifold   $N^{k}\subset M^n$ of dimension  $k$ without boundary  is  {\it locally flat in a point $x\in N^k$} if there exists a neighborhood 
$U(x)\subset M^{n}$ of the point    $x$ and a homeomorphism  $\varphi:U(x)\to
\mathbb{R}^{n}$ such that $\varphi(N^{k}\cap U(x))=
\mathbb{R}^{k}$, where  $\mathbb{R}^k=\{(x_1,...,x_n)\in \mathbb{R}^n|\  x_{k+1}=x_{k+2}=...=x_n=0\}$. 

A manifold   $N^{k}$ is  {\it locally flat} in $M^n$ or 
{\it the submanifold} of the manifold $M^{n}$ if it is locally flat at each  its point.

 If the condition of local flatness fails in a point  $x\in N^k$ then the manifold  $N^k$ is called  {\it wild} and the point $x$ is called {\it the point of wildness}.


A topological space   $X$ is called  {\it $m$-connected} (for $m>0$) if it is non-empty, path-connected and its first $m$  homotopy groups $\pi_i(X)$, $i\in \{1,\dots,m\}$ are trivial. The requirements of being non-empty and path-connected can be interpreted as (-1)-connected and 0-connected correspondingly. 

 A topological space $P$ generated by points of a simplicial complex $K$ with the topology induced from  $\mathbb{R}^n$ is called {\it the polyhedron}. The complex    $K$ is called {\it the partition} or {\it the triangulation} of the polyhedron $P$. 

 A map  $h:P\to Q$  of polyhedra  is called  {\it  piece-vise linear} if there exists partitions  $K, L$ of polyhedra  $P,Q$  correspondingly such that  $h$ move each simplex of the complex $K$ into a simplex of the complex $L$ (see for example  \cite{RuSa}).

A polyhedron $P$ is called  {\it the  piece-vise linear manifold} of dimension  $n$ with boundary if it is a topological manifold with boundary and  for any point  $x\in int\,P$ ($y\in \partial P$)  there is a neighborhood $U_x$ ($U_y$) and a piece-vise linear  homeomorphism  $h_x:U_x\to \mathbb{R}^n$ ($h_y:U_y\to \mathbb{R}^n_+=\{(x_1,...,x_n)\subset \mathbb{R}^n|\ x_1\geq 0\}$).

The following important statement  follows from  Theorem~4 of~\cite{Hu}. 

\begin{st}\label{Huds} Suppose that $N^k, M^n$ are compact  piece-vise linear manifolds of dimension  $k, n$ correspondingly,   $N^k$ is the manifold without boundary,   $M^n$ possibly has a non-empty boundary, $\tilde{e}, e:N^k\to int~M^n$ are homotopic piece-vise linear embeddings, and the following conditions hold:
\begin{enumerate}
\item $n-k\geq 3$;
\item $N^k$ is $(2k-n+1)$-connected;
\item $M^n$ is $(2k-n+2)$-connected.
\end{enumerate}

Then there exists a family of piece-vise linear homeomorphisms   $h_t : M^n\to M^n$, $t\in [0,1]$,  such that   $h_0=id$, $h_1\tilde{e}=e$, $h_t|_{_{\partial M^n}}=id$ for any  $t\in [0,1]$. 
\end{st}

We will say that    a topological submanifold   $N^k\subset M^n$ of the manifold    $M^n$ is  {\it an   essential} if  a homomorphism  ${e_\gamma}_*:\pi_1(N^k)\to \pi_1(M^n)$ induced by an embedding  $e_{_{N^k}}:N^k\to {M}^n$ is the isomorphism.  We will call    an essential manifold  $\beta$ homeomorphic to the circle  $\mathbb{S}^1$ {\it the essential knot}.  

Let  $\beta\in M^n$ be an essential knot and   $h:\mathbb{B}^{n-1}\times \mathbb{S}^1\to {M}^n$ be a topological embedding such that  $h(\{O\}\times \mathbb{S}^1)=\beta$. Call the  image $N_\beta=h(\mathbb{B}^{n-1}\times \mathbb{S}^1)$ {\it the tubular neighborhood} of the knot $\beta$.

\begin{Prop}\label{triv-000} Suppose that $\mathbb{P}^{n-1}$ is either $\mathbb{S}^{n-1}$ or  $\mathbb{B}^{n-1}$,  $\beta_1,...,\beta_k \subset int\,\mathbb{P}^{n-1}\times \mathbb{S}^1$ are essential knots and   $x_1,...,x_k\subset int~\mathbb{P}^{n-1}$ are arbitrary points.  Then there is a homeomorphism   $h: \mathbb{P}^{n-1}\times \mathbb{S}^1 \to \mathbb{P}^{n-1}\times \mathbb{S}^1$ such that     $h(\bigcup\limits_{i=1}^k\beta_i)=\bigcup\limits_{i=1}^k\{x_i\}\times \mathbb{S}^1$ and   $h|_{_{\partial\,\mathbb{P}^{n-1}\times \mathbb{S}^1}}=id$. 
\end{Prop}
\begin{demo}  Put  $b_i=\{x_i\}\times \mathbb{S}^1$, $i\in \{1,...,k\}$. Choose pair-vise disjoint neighborhoods $ U_1,\dots,U_k$  of knots   $\beta_1,\dots,\beta_k$ in $int\,\mathbb{P}^{n-1}\times \mathbb{S}^1$. It follows from  Theorem 1.1 of the paper~\cite{Gluck} that there exists a homeomorphism  $g:\mathbb{P}^{n-1}\times \mathbb{S}^1\to \mathbb{P}^{n-1}\times \mathbb{S}^1$ that is identity outside the set   $\bigcup\limits_{i=1}^{k}U_i$ and such that for any      $i\in \{1,...,k\}$  the set   $g(\beta_i)$ is a subpolyhedron. 

By assumption,    piece-vise linear embeddings  $\tilde{e}:\mathbb{S}^1\times \mathbb{Z}_k\to \mathbb{P}^{n-1}\times \mathbb{S}^1$, $e: \mathbb{S}^1\times \mathbb{Z}_k\to \mathbb{P}^{n-1}\times \mathbb{S}^1$  such that   $\tilde{e}(\mathbb{S}^1\times \mathbb{Z}_k)=\bigcup \limits_{i=1}^kg(\beta_i)$, ${e}(\mathbb{S}^1\times \mathbb{Z}_k)=\bigcup \limits_{i=1}^k b_i$ are homotopic.
By   Statement~\ref{Huds}, there exists a family of piece-vise linear homeomorphisms   $h_t : \mathbb{P}^{n-1}\times \mathbb{S}^1\to  \mathbb{P}^{n-1}\times \mathbb{S}^1$, $t\in [0,1]$,  such that   $h_0=id$, $h_1\tilde{e}=e$, $h_t|_{_{\partial\,\mathbb{P}^{n-1}\times \mathbb{S}^1}}=id$ for any  $t\in [0,1]$. Then  $h_1$ is the desired homeomorphism.  
\end{demo}

The following Statement~\ref{prod}  is proved in the paper~\cite{GrGuMe10} ( see Lemma~2.1). 

\begin{st}\label{prod} Let $h:\mathbb{B}^{n-1}\times \mathbb{S}^1\to  int~\mathbb{B}^{n-1}\times \mathbb{S}^1$ be a topological embedding  such that  $h(\{O\}\times \mathbb{S}^1)=\{O\}\times \mathbb{S}^1$.  Then a manifold 
 $\mathbb{B}^{n-1}\times \mathbb{S}^1\setminus int~h(\mathbb{B}^{n-1}\times \mathbb{S}^1)$ is homeomorphic to the direct product  $\mathbb{S}^{n-2}\times \mathbb{S}^1\times [0,1]$.
\end{st}

\begin{Prop}\label{plus-prod} Suppose that  $Y$ is a topological manifold with boundary, $X$ is a closed component of its boundary, $Y_1$ is a manifold homeomorphic to $X\times [0,1]$, and   $Y\cap Y_1=X$.  Then a manifold $Y\cup Y_1$ is homeomorphic to $Y$. Moreover,   if the manifold $Y$ is homeomorphic to the direct product  $X\times [0,1]$ then there exists a homeomorphism  $h:   X\times [0,1]\to Y\cup Y_1$  such that   $h(X\times \{\frac{1}{2}\})=X$.
 \end{Prop}
\begin{demo} By~\cite{Br62} (Theorem 2),  there exists a topological embedding   $h_0:X\times [0,1]\to Y$ such that  $h_0(X\times \{1\})=X$.  Put $Y_0=h_0(X\times [0,1])$. Let  $h_1:X\times [0,1]\to Y_1$  be a homeomorphism such that $h_1(X\times \{0\})=X=h_0(X\times \{1\})$. 

Define  homeomorphisms  $g:X\times [0,1]\to X\times [0,1]$ , $\tilde{h}_1:X\times [0,1]\to Y_1$, $h:X\times [0,1]\to Y_0\cup Y_1$ by    
$g(x,t)=(h_1^{-1}(h_0(x,1)),t)$, $\tilde{h}_1=h_1g$, 
 
$$h(x,t)=\begin{cases}h_0(x,2t),~
t\in[0,\frac12];
\cr \tilde{h}_1(x,2t-1),~t\in(\frac12;1],\end{cases}$$  and define a homeomoprhism  $H:Y\cup Y_1\to Y$  by  

$$H(x)=\begin{cases}h_0(h^{-1}(x)),~
x\in Y_0\cup Y_1;
\cr x,~x\in Y\setminus Y_0.\end{cases}$$

{To prove the second item of the statement it is enough to put $Y=Y_0$. Then  the homeomorphism   $h:X\times [0,1]\to Y\cup Y_1$ defined above  is the desired one.}
\end{demo}

\begin{Prop}\label{nb-equ} Suppose that $\mathbb{P}^{n-1}$ is either the ball  $\mathbb{B}^{n-1}$ or the sphere  $\mathbb{S}^{n-1}$,  $\beta_1,...,\beta_k\subset int\,\mathbb{P}^{n-1}\times \mathbb{S}^1$  are essential knots,   ${N}_{\beta_1}, ...,N_{\beta_k}\subset \mathbb{P}^{n-1}\times \mathbb{S}^1$ are their  pair-vise disjoint neighborhoods,  $D^{n-1}_1,..., D^{n-1}_k\subset \mathbb{P}^{n-1}$ are  pair-vise disjoint disks, and  $x_1,...,x_k$ are inner points of the disks  $D^{n-1}_1,..., D^{n-1}_k$  correspondingly.  Then there exist a homeomorphism  $h:\mathbb{P}^{n-1}\times \mathbb{S}^1\to \mathbb{P}^{n-1}\times \mathbb{S}^1$  such that $h(\beta_i)=\{x_i\}\times \mathbb{S}^1, h({N}_{\beta_i})=D^{n-1}_i\times \mathbb{S}^1, i\in \{1,...,k\}$ and  $h|_{_{\partial\,\mathbb{P}^{n-1}\times \mathbb{S}^1}}=id$.
\end{Prop}
\begin{demo}  By  Proposition~\ref{triv-000}, there exists a homeomorphism   $h_0: \mathbb{P}^{n-1}\times \mathbb{S}^1\to \mathbb{P}^{n-1}\times \mathbb{S}^1$ such that  $h_0(\beta_i)=\{x_i\}\times \mathbb{S}^1$, $h_0|_{_{\partial\, \mathbb{P}^{n-1}\times \mathbb{S}^1}}=id$. Put  $\tilde{N}_i=h_0(N_{\beta_i})$. By~\cite{Br62},  there exist topological embeddings  $e_i:\mathbb{S}^{n-2}\times\mathbb{S}^1\times [0,1]\to int\,\mathbb{P}^{n-1}\times \mathbb{S}^1$ such that  $e_i(\mathbb{S}^{n-2}\times\mathbb{S}^1\times \{1\})=\partial \tilde{N}_{\beta_i}$, $e_i(\mathbb{S}^{n-2}\times\mathbb{S}^1\times [0,1])\cap e_j(\mathbb{S}^{n-2}\times\mathbb{S}^1\times [0,1])=\emptyset$ for $i\neq j, i,j\in \{1,...,k\}$. Put $U_i=e_i(\mathbb{S}^{n-2}\times\mathbb{S}^1\times [0,1])\cup \tilde{N}_i$. 

 Suppose that $D^{n-1}_{0,1},...,D^{n-1}_{0,k}, {D}^{n-1}_{1,1},...,{D}^{n-1}_{1,k}\subset \mathbb{P}^{n-1}$ are disks such that   $x_i\subset int\,{D}^{n-1}_{j,i}$, $D^{n-1}_{j,i}\subset int\, D^{n-1}_i$, $j\in \{0,1\}$, ${D}^{n-1}_{0,i}\subset int\,D^{n-1}_{1,i}$, and $D^{n-1}_{1,i}\times \mathbb{S}^1\subset int\,\tilde{N}_i$. 

By  Proposition~\ref{prod}, every set  $\tilde{N}_i\setminus (int\,D^{n-1}_{1,i}\times\mathbb{S}^1)$, $(D^{n-1}_{1,i}\setminus int\,D^{n-1}_{0,1})\times \mathbb{S}^1$ is homeomorphic to the direct product $\mathbb{S}^{n-2}\times \mathbb{S}^1\times [0,1]$. By  Proposition~\ref{plus-prod}, there exists a homeomorphism  $g_i:\mathbb{S}^{n-2}\times \mathbb{S}^1\times [0,1]\to U_i\setminus int\,D^{n-1}_{0,i}\times \mathbb{S}^1$ such that $g_i(\mathbb{S}^{n-2}\times \mathbb{S}^1\times \{t_1\})=\partial\, \tilde{N}_i$, $g_i(\mathbb{S}^{n-2}\times \mathbb{S}^1\times \{t_2\})=\partial\, D^{n-1}_{1,i}\times \mathbb{S}^1$ for some  $t_1,t_2\subset (0,1)$.   Let $\xi:[0,1]\to [0,1]$ be a homeomorphism that is identity on the ends of the interval  $[0,1]$ and such that $\xi(t_1)=t_2$.  Define a homeomorphism  $\tilde{g}_i:\mathbb{S}^{n-2}\times \mathbb{S}^1\times [0,1]\to  \mathbb{S}^{n-2}\times \mathbb{S}^1\times [0,1]$  by  $\tilde{g}_i(x,t)=(x,\xi(t))$. 

Define a homeomorphism $h_i:\mathbb{P}^{n-1}\times \mathbb{S}^1\to \mathbb{P}^{n-1}\times \mathbb{S}^1$  by

$$h_i(x)=\begin{cases}
g_i(\tilde{g}_i(g^{-1}_i(x))),~x\in U_i\setminus int\,D^{n-1}_{0,i}\times \mathbb{S}^1;
\cr x,~x\in (\mathbb{P}^{n-1}\times \mathbb{S}^1\setminus U_i).\end{cases}$$

The superposition  $\eta=h_k\cdots h_1h_0$ maps every knot  $\beta_i$ into  the knot  $\{x_i\}\times \mathbb{S}^1$,   the  neighborhood $N_{\beta_i}$ into the set    $D^{n-1}_{1,i}\times \mathbb{S}^1$, and keeps  the set  $\partial\, \mathbb{P}^{n-1}\times \mathbb{S}^1$ fixed.  Construct a homeomorphism  $\Theta:\mathbb{P}^{n-1}\times \mathbb{S}^1\to \mathbb{P}^{n-1}\times \mathbb{S}^1$ that  be identity on the set  $\partial\,\mathbb{P}^{n-1}\times \mathbb{S}^1$ and on the knots  $\{x_1\}\times \mathbb{S}^1,...,\{x_k\}\times \mathbb{S}^1$ and  move the set  $D^{n-1}_{1,i}\times \mathbb{S}^1$ into  the set  $D^{n-1}_{i}\times \mathbb{S}^1$ for every  $i\in\{1,...,k\}$. It follows from the Annulus Theorem\footnote{The Annulus Theorem  states that the closure of an open domain on the sphere $S^{n+1}$ bounded by two disjoint locally flat spheres $S^{n}_1, S^{n}_2$ is homeomorphic to the annulus $\mathbb{S}^{n}\times [0,1]$. In dimension    2 it was proved by Rado in  1924,  in dimension 3 --- by Moise in 1952, in dimension   4 --- by Quinn in 1982, and in dimension 5 and greater --- by  Kirby in 1969.} that the set  $D^{n-1}_{i}\setminus int\,D^{n-1}_{1,i}$ is homeomorphic to the annulus  $\mathbb{S}^{n-2}\times [0,1]$. Then apply the construction similar to one  described above  to define  a homeomorphism  $\theta:\mathbb{P}^{n-1}\to \mathbb{P}^{n-1}$ such that  $\theta (x_i)=x_i$, $\theta(D^{n-1}_{i})=D^{n-1}_{1,i}$, $\theta|_{_{\partial\,\mathbb{P}^{n-1}}}=id$. Put   $\Theta(x,t)=(\theta^{-1}(x),t)$, $x\in \mathbb{P}^{n-1}, t\in \mathbb{S}^1$. Then   $h=\Theta \eta$ is the desired homeomorphism.   
\end{demo}
\begin{Cor}\label{blue}
 If $N\subset \mathbb{S}^{n-1}\times \mathbb{S}^1$ is a tubular neighborhood of an essential knot than the manifold  $(\mathbb{S}^{n-1}\times \mathbb{S}^1)\setminus int~N$ is homeomorphic to the direct product  $\mathbb{B}^{n-1}\times \mathbb{S}^1$.
\end{Cor}

\subsection{A surgery of the manifold  $\mathbb{S}^{n-1}\times \mathbb{S}^1$ along an essential  submanifold homeomorphic to  $\mathbb{S}^{n-2}\times \mathbb{S}^1$}\label{surg}

Recall that we put  $\mathbb{Q}^n=\mathbb{S}^{n-1}\times \mathbb{S}^1$. Suppose that  $N\subset \mathbb{Q}^{n}$  is an essential  submanifold homeomorphic to  $\mathbb{B}^{n-1}\times \mathbb{S}^1$,  $T=\partial N$, and  $e_T:\mathbb{S}^{n-2}\times \mathbb{S}^1\times [-1;1]\to \mathbb{Q}^{n}$ is a topological embedding such that  $e_T(\mathbb{S}^{n-2}\times \mathbb{S}^1\times \{0\})=T$. Put $K=e_T(\mathbb{S}^{n-2}\times \mathbb{S}^1\times [-1;1])$ and denote by  $N_+, N_-$ connected components of the set  $\mathbb{Q}^{n}\setminus int\, K$. It follows from  Propositions~\ref{nb-equ},~\ref{plus-prod} that the manifolds  $N_+, N_-$ are homeomorphic to  $\mathbb{B}^{n-1}\times \mathbb{S}^1$. Let  $N'_+, N'_-$   manifolds homeomorphic to  $\mathbb{B}^{n-1}\times \mathbb{S}^1$.  Denote by $\psi_\delta: \partial\,N_{\delta}\to \partial\, N'_{\delta}$ an arbitrary homeomorphism reversing the natural orientation, by $Q_\delta$  a manifold obtained by gluing the      manifolds  $N_\delta$ and $N'_{\delta}$ by means of homeomorphism $\psi_{\delta}$, and by     $\pi_{{\delta}}:(N_\delta \cup N'_\delta)   \to Q_{\delta}$ the natural projection,  $\delta\in \{+,-\}$. 

We will say that the   manifolds  $Q_+, Q_-$ are obtained from  $\mathbb{Q}^{n}$ {\it by the surgery along the submanifold  $T$}.

Note that $\mathbb{S}^{n-2}\times \mathbb{S}^1$ is the boundary of   $\mathbb{B}^{n-1}\times \mathbb{S}^1$. 
By~\cite{Max} (Theorem 2), the following statement holds.

\begin{st}\label{maxx} Let $\psi:\mathbb{S}^{n-2}\times \mathbb{S}^1\to \mathbb{S}^{n-2}\times \mathbb{S}^1$ be an arbitrary homeomorphism. Then there exists a homeomorphism  $\Psi:\mathbb{B}^{n-1}\times \mathbb{S}^1\to \mathbb{B}^{n-1}\times \mathbb{S}^1$ such that  $\Psi|_{_{\mathbb{S}^{n-2}\times \mathbb{S}^1}}=\psi|_{_{\mathbb{S}^{n-2}\times \mathbb{S}^1}}$.
\end{st}

\begin{Prop}\label{forward} The manifolds $Q_+$, $Q_-$ are homeomorphic to   $\mathbb{Q}^n$. 
\end{Prop}
\begin{demo}  Let   $D^{n-1}\subset \mathbb{S}^{n-1}$ be an arbitrary disk, $\mathbb{N}_\delta=D^{n-1}\times \mathbb{S}^1$ and  $h_\delta: \pi_{\delta}(N_\delta)\to \mathbb{N}_\delta$ be an arbitrary homeomorphism.  Put $\tilde{\psi}_\delta=h_\delta \pi_{\delta} \psi_\delta\pi^{-1}_{\delta}h^{-1}_\delta|_{\partial\,\mathbb{N}_\delta}$. Due to  Proposition~\ref{maxx} a homeomorphism  $\tilde{\psi}_\delta$ can  extend up to  a homeomorphism  $h'_\delta:\pi_\delta(N'_\delta)\to \mathbb{Q}^n\setminus int\,\mathbb{N}_\delta$. Then a map $H_\delta:Q_\delta\to \mathbb{Q}^n$ defined by   $H_\delta(x)=h_\delta(x)$ whenever $x\in \pi_\delta({N_\delta})$ and  $H_\delta(x)=h'_\delta(x)$ whenever  $x\in \pi_\delta ({N}'_\delta)$ is the desired homeomorphism.
\end{demo} 

\subsection{A surgery of manifolds homeomorphic to $\mathbb{S}^{n-1}\times \mathbb{S}^1$ along essential knots}
\label{inv-surg}

 Let   $Q^n_1,\dots,Q^{n}_{k+1}$ be manifolds homeomorphic to $\mathbb{Q}^n$ . Denote  by  $\beta_1,...,\beta_{2k}\subset \bigcup \limits_{i=1}^{k+1}Q^n_{i}$  essential knots such that for any $j\in \{1,...,k\}$ knots $\beta_{2j-1},\beta_{2j}$ belongs to distinct manifolds from the union  $\bigcup \limits_{i=1}^{k+1}Q^n_{i}$ and every manifold  $Q^n_i$ contains at least one knot from the set  $\beta_1,...,\beta_{2k}$. Let   $N_{\beta_1},...,N_{\beta_{2k}}$ be  tubular neighborhoods of the knots  $\beta_1,...,\beta_{2k}$ correspondingly.  

Let $K_1,...,K_k$ be   manifolds homeomorphic to the direct product  $\mathbb{S}^{n-2}\times \mathbb{S}^1\times [-1;1]$. For every  $j\in \{1,\dots,k\}$ denote by   $T_j\subset K_j$ a manifold homeomorphic to  $\mathbb{S}^{n-2}\times \mathbb{S}^1$ that cuts  $K_j$ into two connected components whose closures are homeomorphic to  $\mathbb{S}^{n-2}\times \mathbb{S}^1\times [0;1]$,  and by     $\psi_{j}:\partial N_{2j-1}\cup \partial  N_{2j} \to \partial K_j$ an arbitrary reversing the natural orientation homeomorphism. 

Glue manifolds  $\widetilde{Q}=\bigcup \limits_{i=1}^{k+1}Q^n_{i}\setminus~\bigcup\limits_{\nu=1}^{2k}int~N_\nu$ and  $K=\bigcup \limits_{j=1}^{k}K_j$ by means of the homeomorphisms  $\psi_{1},...,\psi_{k}$,   denote by $Q$ the obtained manifold and by  $\pi:\widetilde{Q} \cup K\to Q$ the  natural projection. We will say that the manifold $Q$ is obtained from   $Q^n_{1},...,Q^n_{k+1}$ by  {\it the surgery along knots  $\beta_1,...,\beta_{2k}$} and call every pair  $\beta_{2j-1}, \beta_{2j}$ {\it the binding pair}, $j\in \{1,2,...,k\}$.

\begin{Prop}\label{element}  The manifold  $Q$ is  homeomorphic to $\mathbb{Q}^{n}$ and every manifold  $\pi(T_j)$ cuts  $Q$ into two connected components whose closures are homeomorphic to  $\mathbb{B}^{n-1}\times \mathbb{S}^1$.
\end{Prop}
\begin{demo} Prove the proposition by induction on  $k$. Consider the case  $k=1$.  Due to Propositions~\ref{nb-equ},~\ref{plus-prod} manifolds  $\widetilde{N}_1=Q^n_1\setminus int\,N_1$, $\widetilde{N}_2=Q^n_2\setminus int\,N_2$, $\widetilde{N}_1\bigcup_{\psi_1|_{_{\partial\,N_1}}}K_1$  are homeomorphic to the direct product  $\mathbb{B}^{n-1}\times \mathbb{S}^1$.  By definition, the manifold  $T_1$ cuts the manifold  $K_1$ into two connected components whose closures are homeomorphic to  $\mathbb{Q}^{n-1}\times [0,1]$. It follows  from  Proposition~\ref{plus-prod} that $T_1$ cuts   $\widetilde{N}_1\bigcup_{\psi_1|_{_{\partial\,N_1}}}K_1$ into two connected components such that the closure of one of which, denote it by $N$,   is homeomorphic to  $\mathbb{B}^{n-1}\times \mathbb{S}^1$ and the closure of another is homeomorphic to  $\mathbb{Q}^{n-1}\times [0,1]$.   Suppose that  $D^{n-1}_0\subset \mathbb{S}^{n-1}$ is an arbitrary disk,  $N_0=D^{n-1}_0\times \mathbb{S}^1$ and   $h_0: \pi(\widetilde{N}_1\bigcup K_1)\to N_0$  is an arbitrary homeomorphism.  Put  $\tilde{\psi}_1=h_0 \pi \psi^{-1}_1\pi^{-1}h_0^{-1}|_{\partial\,N_0}$. In virtue of Proposition~\ref{maxx} a homeomorphism  $\tilde{\psi}$ can be extended up to  a homeomorphism  $h_1:\pi(\widetilde{N}_2)\to \mathbb{Q}^n\setminus int\,N_0$. Then the map  $h:Q\to \mathbb{Q}^n$ defined by   $h(x)=h_0(x)$ for $x\in \pi(\widetilde{N}_1\bigcup K_1)$ and $h(x)=h_1(x)$ for $x\in \pi (\widetilde{N}_2)$ is the desired homeomorphism.  The manifold $\pi(T_1)$ cuts  $Q$ into  two connected components such that  the closure of one of them is   $\pi (N)$ which  is  homeomorphic to  $\mathbb{B}^{n-1}\times \mathbb{S}^1$. By  Corollary~\ref{blue}, the closure of another connected component is also homeomorphic to $\mathbb{B}^{n-1}\times \mathbb{S}^1$.
  
Suppose that the statement is true for all $\lambda=k$ and show that it is true also for  $\lambda=k+1$. Since  $2k\geq k+1$ we have  that  there exists at least one manifold  among the manifolds  $Q^n_1,...,Q^n_{\lambda+1}$, say  $Q^n_{\lambda+1}$, containing exactly one knot  from the set  $\beta_1,...,\beta_{2k}$  (if every of that manifolds would contain no less than two knots, then the total  number of all  knots be no less than  $2k+2$). Let $\beta_{2\lambda}\subset Q^n_{\lambda+1}$, $\beta_{2\lambda-1}\subset Q^n_i$, $i\in \{1,\dots,\lambda\}$, be a binding pair. By  the induction hypothesis and  Corollary~\ref{blue}, the manifold $Q_{\lambda}$ obtained by the surgery of manifolds $Q^n_1,...,Q^n_{\lambda}$ along knots $\beta_1,...,\beta_{2\lambda-2}$ is homeomorphic to  $\mathbb{Q}^n$; the projection of every manifold  $(T_j)$ cuts $Q_{\lambda}$ into two connected components such that   the closure  of each of which is homeomorphic to  $\mathbb{B}^{n-1}\times \mathbb{S}^1$; and the projection of the knot  $\beta_{2\lambda-1}$ is the essential knot.   Now apply the surgery to manifolds  $Q_{\lambda}$, $Q^n_{\lambda+1}$ along knots  $\pi(\beta_{2\lambda-1})$, $\beta_{2\lambda}$ and use the first step  arguments to obtain the desired statement. 
 \end{demo}


\subsection{Proof of Lemma~\ref{simple}} 

{\bf Step  1.} {\it  Proof of the fact  that the manifold  $\widehat{V}_f$  is homeomorphic to  $\mathbb{Q}^n$ and every connected component  $\mathcal{Q}^{n-1}$ of the set  $\hat{L}^u_f\cup \hat{L}^s_f$ cuts  $\widehat{V}_f$  into two connected components whose  closures are homeomorphic to  $\mathbb{B}^{n-1}\times \mathbb{S}^1$.}

Put $k_i=|\Omega^i_f|$, $i\in \{0,1,n-1,n\}$. Due to  Statement~\ref{smale} and the fact that the closure of every separatrix of dimension $(n-1)$ cuts the ambient sphere  $S^n$ into two connected components one gets $k_0=k_1+1$, $k_n=k_{n-1}+1$. 

 Denote by   $\beta_1,...,\beta_{2k_1}$  the essential knots in the set $\widehat{V}=\bigcup \limits_{\omega\in \Omega^0_f}\widehat{V}^s_\omega$ which  are projections  (by means of  $p_{\widehat{V}}$) of all one-dimension unstable separatrices   of the diffeomorphism $f$. Without loss of generality assume that  knots  $\beta_{2j-1},\beta_{2j}$ are the projection of the separatrices of the same saddle point  $\sigma_j\in \Omega^1_f$, $j\in \{1,...,k_1\}$. 

It follows from  Statement~\ref{smale} that every manifold  $\widehat{V}^s_\omega$ contains at least one knot from the set  $\beta_1,...,\beta_{2k_1}$.   Since stable and unstable manifolds of different saddle points do not intersect we have that  for any $j\in \{1,...,k_1\}$ knots $\beta_{2j-1},\beta_{2j}$ belong to distinct connected components of   $\widehat{V}$.  Indeed, if one suppose that  $\beta_{2j-1},\beta_{2j}\subset \widehat{V}^s_\omega $  for some $j, \omega$, then the set  $cl\,W^u_{\sigma_j}=W^u_{\sigma_j}\cup \omega$ is homeomorphic to the circle.  Since $cl\, W^s_{\sigma_j}$ divides the sphere $S^n$ into two parts and intersect the circle $cl\,W^u_{\sigma_j}$ at the  point $\sigma_j$ we have that there exists  at least one  point in $cl\, W^s_{\sigma_j}\cap cl\,W^u_{\sigma_j}$ different from  $\sigma_j$. This fact contradicts to the item 1 of Statement~\ref{homo}.

Let  $N_{\sigma_j}$, $\chi_{\sigma_j}: N_{\sigma_j}\to \mathbb{U}$ be  the neighborhood  of the point  $\sigma_j$ and the homeomorphism defined in  Proposition~\ref{grgupo-ado}. Further we use denotations of the sections~\ref{MS-properties}, \ref{aux}. Denote by  $N_{2j-1}, N_{2j}$  the connected components of the set $\widehat{N}^{u}_{\sigma_j}$ containing knots  $\beta_{2j-1}, \beta_{2j}$ correspondingly. Let  $\psi:\partial\widehat{\mathbb{N}}^{u}\to \partial \widehat{\mathbb{N}}^{s}$ be a homeomorphism such that $\psi p_u|_{\partial \mathbb{U}}=p_s|_{\partial \mathbb{U}}$. Put $K_j= \widehat{N}^{s}_{\sigma_j}$,  $T_j=\widehat{V}^s_{\sigma_j}$
and define homeomorphisms  $\varphi_{u,j}:N_{2j-1}\cup  N_{2j}\to \widehat{\mathbb{N}}^u$,  $\varphi_{s,j}: K_j\to \widehat{\mathbb{N}}^s$, $\psi_{j}:\partial N_{2j-1}\cup \partial  N_{2j} \to \partial K_j$ by  

$$ \varphi_{u,j}=p_u\chi_{\sigma_j}p^{-1}_{\widehat{V}_f}|_{N_{2j-1}\cup N_{2j}},$$

 $$\varphi_{s,j}=p_s\chi_{\sigma_j}p^{-1}_{\widehat{V}_f}|_{K_j}, $$

$$\psi_{j}=\varphi^{-1}_{s,j}\psi\varphi_{u,j}|_{\partial N_{2j-1}\cup \partial N_{2j}}, $$

  and denote by $$\Psi:\bigcup\limits_{j=1}^{k_1}(\partial N_{2j-1}\cup \partial N_{2j})\to \bigcup\limits_{j=1}^{k_1}K_j$$ the homeomorphism such that $$\Psi|_{_{\partial N_{2j-1}\cup \partial N_{2j}}}=\psi_j|_{_{\partial N_{2j-1}\cup \partial N_{2j}}}.$$

	Since $$V_f=\left(\bigcup\limits_{\omega\in \Omega^0_f}V^s_\omega\setminus \left(\bigcup\limits_{\sigma\in \Omega^1_f}V^u_\sigma\right)\right)\bigcup \left(\bigcup\limits_{\sigma\in \Omega^1_f}V^s_\sigma\right) =\left(V_f\setminus \left(\bigcup\limits_{\sigma\in \Omega^1_f}N^u_\sigma\right)\right)\bigcup \left(\bigcup\limits_{\sigma\in \Omega^1_f}N^s_\sigma\right)$$ it follows that

$$\widehat{V}_f=\left(\widehat{V}_f\setminus \left(\bigcup\limits_{\sigma\in \Omega^1_f}\widehat{N}^u_\sigma\right)\right)\cup _{_\Psi}\left(\bigcup\limits_{\sigma\in \Omega^1_f}\widehat{N}^s_\sigma\right)=\left(\widehat{V}_f\setminus \left(\bigcup\limits_{j=1}^{2k_1}{N}_{j}\right)\right)\cup_{_\Psi} \left(\bigcup\limits_{j=1}^{k_1}K_j\right).$$

 So, the manifold  $\widehat{V}_f$  is obtained from  $\bigcup \limits_{\omega\in \Omega^0_f}\widehat{V}^s_\omega$   by the surgery along knots  $\beta_{1},...,\beta_{2k_1}$. Due to  Proposition~\ref{element}, the manifold  $\widehat{V}_f$ is homeomorphic to   $\mathbb{Q}^n$ and every connected component of the set   $\hat{L}^s_f$  cuts the set  $\widehat{V}_f$ into two connected components such that the closure of each of which is homeomorphic to   $\mathbb{B}^{n-1}\times \mathbb{S}^1$. 
		
	From the other hand $$V_f=\left(\bigcup\limits_{\alpha\in \Omega^n_f}V^u_\alpha \setminus \left(\bigcup\limits_{\sigma\in \Omega^{n-1}_f}V^s_\sigma\right)\right)\bigcup \left(\bigcup\limits_{\sigma\in \Omega^{n-1}_f}V^u_\sigma\right) =\left( V_f\setminus \left(\bigcup\limits_{\sigma\in \Omega^{n-1}_f}N^s_\sigma\right)\right)\bigcup \left(\bigcup\limits_{\sigma\in \Omega^{n-1}_f}N^u_\sigma\right).$$ Similar to previous arguments one can conclude that  the set  $\widehat{V}_f$ is obtained from  $\bigcup \limits_{\alpha\in \Omega^n_f}\widehat{V}^u_\alpha$   by the surgery along the projections of all  one-dimensional stable separatrices of the saddle points of the diffeomorphism  $f$. In virtue of  Proposition~\ref{element} every   connected component of the set  $\hat{L}^u_f$  cuts the set  $\widehat{V}_f$ into two connected components such that the closure of each of which is homeomorphic to $\mathbb{B}^{n-1}\times \mathbb{S}^1$.

{\bf Step 2.}	{\it Proof of the fact  that there is a set $\widehat{\mathbb{L}}_{m_f}\subset \mathbb Q^n$ and a homeomorphism  $\hat\varphi:\widehat V_f\to\mathbb Q^n$  such that   $\hat\varphi(\hat{L}^s_{f}\cup \hat{L}^u_{f})=\widehat{\mathbb{L}}_{m_f}$.} 
	
	Denote by  $\mathcal{Q}^{n-1}_1,...,\mathcal{Q}^{n-1}_{_{k_1+k_{n-1}}}$ all 
elements of the set $\hat{L}^s_{f}\cup \hat{L}^u_{f}$ and suppose that  $\mathcal{Q}^{n-1}_1$ is  an element such that  all elements of the set   $\hat{L}^s_{f}\cup \hat{L}^u_{f}\setminus \mathcal{Q}^{n-1}_1$ are contained exactly in one of the connected component of the manifold $\widehat{V}_f\setminus \mathcal{Q}^{n-1}_1$. Denote by $N_1$ the closure of this connected component. By  Step 1,  $N_1$ is homeomorphic to  $\mathbb{B}^{n-1}\times \mathbb{S}^1$. By Proposition~\ref{nb-equ},  there exists a disk $D^{n-1}_1\subset \mathbb{S}^{n-1}$ and a  homeomorphism  $\psi_0: \hat{V}_f\to \mathbb{Q}^n$ such that $\psi_0(N_1)=D^{n-1}_1\times \mathbb{S}^1$. If $k_1+k_{n-1}=1$ then the proof is complete and $\hat\varphi=\psi_0$,   $\widehat{\mathbb{L}}_{m_f}=\partial{D}^{n-1}_1\times \mathbb{S}^1$. 
		
		 Let $k_1+k_{n-1}>1$. Denote the images of  $\mathcal{Q}^{n-1}_1,...,\mathcal{Q}^{n-1}_{_{k_1+k_{n-1}}}$  under the homeomorphism  $\psi_0$  by the same symbols as their originals. For $i\in \{2,\dots,k_1+k_{n-1}\}$ denote by  $N_i$ the connected component of the set  $\mathbb{Q}^n\setminus \mathcal{Q}^{n-1}_i$  contained in the set   $D^{n-1}_1\times \mathbb{S}^1$. Without loss of generality suppose that the  numeration of the sets $\mathcal{Q}^{n-1}_1,...,\mathcal{Q}^{n-1}_{_{k_1+k_{n-1}}}$  is chosen in such a way that there exist a number  $l_1\in [2,  k_1+k_{n-1}]$ and  pair-vise disjont sets  $N_2,\dots,N_{l_1}$ such that     $\bigcup\limits_{i=2}^{l_1}N_i=\bigcup\limits_{i=2}^{k_1+k_{n-1}}N_i$.  Choose in the interior of the disk   $D^{n-1}_1$  arbitrary pair-vise disjoint disks  $D^{n-1}_2,\dots,D^{n-2}_{l_1}$.  Due to Proposition~\ref{nb-equ} there exists a homeomorphism $\psi_1:\mathbb{Q}^n\to \mathbb{Q}^n$ such that  $\psi_1|_{\mathbb{Q}^n\setminus int\,D^{n-1}_1\times \mathbb{S}^1}=id$, $\psi_1(N_i)=D^{n-1}_i\times \mathbb{S}^1$, $i\in \{2,\dots,l_1\}$. If   $l_1=k_1+k_{n-1}$ then the proof is complete and  $\hat\varphi=\psi_1\psi_0$,  $\widehat{\mathbb{L}}_{m_f}=\bigcup\limits_{i=1}^{l_1}\partial\,D^{n-1}_i\times \mathbb{S}^1$.
		
		Let  $l_1< k_1+k_{n-1}$. Denote the images of $\mathcal{Q}^{n-1}_1,\dots,\mathcal{Q}^{n-1}_{_{k_1+k_{n-1}}}$ and  $N_1,\dots,N_{k_1+k_{n-1}}$   under  the homeomorphism  $\psi_1$ by the same symbols as their originals. Put $\mathcal{N}=\bigcup\limits_{i=l_1+1}^{k_1+k_{n-1}}N_{i}$. 

 If  for fixed  $i\in \{2,...,l_1\}$ the set   $N_i$  has  non-empty intersection with  the set $\mathcal{N}$, then   denote by  $l_i, \tilde{k}_i$,   $l_i\leq \tilde{k}_i$, the   positive numbers  such that    $N_{i,1},...,N_{i,\tilde{k}_i}$ are   all elements from  $N_i\cap \mathcal{N}$ and  $N_{i,1},...,N_{i,l_i}$ are pair-vise disjoint  elements from  $N_i\cap \mathcal{N}$  such that 
		$\bigcup\limits_{j=1}^{l_i}N_{i,j}=\bigcup\limits_{j=2}^{\tilde{k}_i}N_{i,j}$. Choose in the interior of the every  disk   $D^{n-1}_i$  pair-vise disjoint disks  $D^{n-1}_{i,1},\dots D^{n-1}_{i,l_i}$.   It follows from  Proposition~\ref{nb-equ}  that there exists a homeomorphism  $\psi_i:\mathbb{Q}^n\to \mathbb{Q}^n$  such that  $\psi_i|_{\mathbb{Q}^n\setminus int\,N_i}=id$, $\psi_i(N_{i,j})=D^{n-1}_{i,j}\times \mathbb{S}^1$, $j\in \{1,\dots,l_i\}$, $i\in \{2,...,l_1\}$. If $N_i\cap \mathcal{N}=\emptyset$, put  $\psi_i=id$.  

If   $l_i=\tilde{k}_i$ for any $i\in\{2,\dots,l_1\}$ such that the numbers $l_i, \tilde{k}_i$ are defined, then  the proof is complete and  $\hat{\varphi}=\psi_{l_1}\psi_{l_1-1}\cdots\psi_1$, $\widehat{\mathbb{L}}_{m_f}= \bigcup\limits_{i=1}^{l_1}\bigcup\limits_{j=1}^{l_i}\partial D^{n-1}_{i,j}\times \mathbb{S}^1$.  Otherwise, continue the process and after finite number of steps get the desired set $\widehat{\mathbb{L}}_{m_f}$ and the desired homeomorphism   $\hat\varphi$ as a superposition of all constructed homeomorphisms.


	\section{Embedding of diffeomorphisms from the class  $G(M^n)$ into topological flows}\label{flow}
	
	\subsection{Free and properly discontinuous action of a group of maps}\label{vsp}
In this section we collect an axillary facts on properties of the transformation group   $\{g^{n}, n\in \mathbb{Z}\}$ which is an  infinite cyclic group acting freely and properly discontinuously  on a topological (in general, non-compact) manifold  $X$ and generated by a homeomorphism  $g:X \to X$\footnote{A group 
$\mathcal G$ {\it acts} on the manifold $X$ if there is  a map 
$\zeta:\mathcal G\times X\to X$ with the following properties:

1) $\zeta(e,x)=x$ for all  $x\in X$, where  $e$
is the identity element of the group  $\mathcal G$;

2) $\zeta(g,\zeta(h,x))=\zeta(gh,x)$ for all  $x\in X$ and $g,h\in
\mathcal G$.

A group  $\mathcal G$  acts {\it freely} on a manifold  $X$ if for any different  $g,h\in X$ and for any point  $x\in X$ an inequality  $\zeta(g,x)\neq \zeta(h,x)$ holds.

A group $\mathcal G$  acts 
{\it properly discontinuously} on the manifold  $X$ if for every compact subset  $K\subset X$ the set of elements  $g\in\mathcal G$ such that 
$\zeta(g,K)\cap K\neq\emptyset$ is finite.}.

Denote by  $X/g$  the orbit space of the action of the group  $\{g^{n}, n\in \mathbb{Z}\}$ and by  $p_{_{X/g}}:X\to {X}/g$ the natural projection. In virtue of~\cite{Th01} (Theorem~3.5.7  and Proposition~3.6.7)  the natural projection    $p_{_{X/g}}\colon
X\to {X}/g$  is a covering map and the space  ${X}/g$ is a manifold.

Denote by  $\eta_{_{X/g}}:\pi_{1}({X}/g)\to \mathbb{Z}$ a homeomorphism defined in the following way. Let $\hat{c}\subset {X}/g$ be a loop non-homotopic to zero in   ${X}/g$ and   $[\hat{c}]\in \pi_1({X}/g)$ be a homotopy class of  $\hat{c}$.  Choose an arbitrary point  $\hat{x}\in \hat{c}$, denote by   $p_{_{X/g}}^{-1}(\hat{x})$  the complete inverse image of  $\hat{x}$, and fix a point  $\tilde{x}\in p_{{_{X/g}}}^{-1}(\hat{x})$. As   $p_{_{X/g}}$ is the covering map then   there is a unique path  $\tilde {c} (t)$  beginning at the point $\tilde {x}$ ($\tilde
{c} (0) = \tilde x$) and covering the loop $c$ (such that  $p_{_{X/g}}(\tilde {c} (t))=
\hat{c}$). Then there  exists the element  $n\in\mathbb {Z}$ such that 
$\tilde{c}(1)=f^n (\tilde {x})$. Put  $\eta_{_{X/g}}([\hat{c}]) =
n$. It follows from~\cite{Ko} (гл.~18)  that the homomorphism  $\eta_{_{X/g}}$ is an epimorphism. 

The next statement~\ref{fact} can be found in~\cite{Ko}  (Theorem~5.5) and~\cite{BoGrPo05} (Propositions 1.2.3 и 1.2.4). 

\begin{st}\label{fact} Suppose that  $X$, $Y$ are connected topological manifolds and  $g:X\to X$, $h:Y\to Y$ are homeomorphisms such that groups 
$\{g^n,~n\in\mathbb {Z} \} $, $ \{h^n,~n\in\mathbb {Z}\} $
acts freely and properly discontinuously on  $X$, $Y$
correspondingly. Then:
\begin{itemize}
\item[1$)$] if  $\varphi:X\to Y$ is a homeomorphism  such that $h=\varphi g \varphi^{-1}$  and $\varphi_*:\pi_1 ({X}/g)\to \pi_1 ({Y}/h)$ is the induced homomorphism, then a map 
$\widehat{\varphi}:{X}/g\to {Y}/h$ defined by 
 $\widehat{\varphi}=p_{{_{Y/h}}}\varphi p^{-1}_{{_{X/g}}}$ is a homeomorphism and  $\eta_{_{X/g}}=\eta_{_{Y/h}}\varphi_*$;

\item[2$)$] if   $\widehat{\varphi}:{X}/g\to {Y}/h$ is a homeomorphism such that 
$\eta_{_{X/g}}=\eta_{_{Y/h}}\varphi_*$ and  
 $\hat{x}\in {X}/g$, $\tilde{x}\in p ^{-1}_{_{X/g}}(x)$, $y=\widehat{\varphi}(x)$,  $\tilde{y}\in p^{-1}{_{Y/h}}(y)$, then there exists a  unique homeomorphism  $\varphi:X\to
Y$ such that   $h=\varphi g \varphi^{-1}$  and   $\varphi (\tilde x) = \tilde y$.
\end{itemize}
\end{st}

\subsection{Proof of  Theorem~\ref{emb-to-flow}}
	
	Suppose that a Morse-Smale  diffeomorphism 	$f:S^n\to S^n$ has no heteroclinic intersection and satisfy Palis conditions. To prove the theorem it is enough to construct  a topological flow  $X^t_f$ such that its time one  map $X^1_f$ belongs to the class  $G(S^n)$ and the scheme  $S_{X^1_f}$ is equivalent to the scheme  $S_f$ (see Section~\ref{sort}).

{\bf Step 1.} It follows from  Lemma~\ref{simple} and  Proposition~\ref{fact} that there exists a homeomorphism $\psi_{_f}:V_f\to \mathbb{S}^{n-1}\times \mathbb{R}$ such that:  

1) $f|_{V_f}=\psi_{_f}^{-1}a\psi_{_f}$, where $a$ is the time one  map of the flow  $a^t(x,s)=(x,s+t)$, $x\in S^{n-1}, s\in \mathbb{R}$;

2) for $(n-1)$-dimensional separatrix $l_\sigma$ of an arbitrary  saddle point $\sigma\in \Omega_f$ there exists a sphere $S^{n-2}_{\sigma}\subset \mathbb{S}^{n-1}$ such that  $\psi_{_f}(l_\sigma)=\bigcup\limits_{t\in \mathbb{R}}a^t(S^{n-2}_{\sigma})$.

Recall that we denote by  $L^s_f$ and  $L^u_f$ the union of all $(n-1)$-dimensional stable and unstable separatrices of the diffeomorphism $f$  correspondingly. Put $\mathbb L^s=\psi_f(L^s_f)$, $\mathbb L^u=\psi_f(L^u_f)$. Then  ${\mathbb L}^\delta$ is the union of pair-vise disjoint  cylinders $\tilde{Q}^\delta_1\cup\dots\cup\tilde{Q}^\delta_{k^\delta},~\delta\in\{s,u\}$. Denote by  $N({\mathbb L}^\delta)=N(\tilde{Q}^\delta_1)\cup\dots\cup N(\tilde{Q}^\delta_{k^\delta})$  the set of their  pair-vise disjoint  closed tubular neighborhoods such that  $N(\tilde{Q}^\delta_i)=K^\delta_i\times\mathbb R$, where $K^\delta_i\subset\mathbb S^{n-1}$ is an annulus of dimension  $(n-1)$,  $i=1,\dots,k^\delta$. 

 Define a flow  $a^t_1$ on the set  $\mathbb{U}=\{(x_{1},...,x_{n})\in
\mathbb{R}^{n}|\ x_{1}^{2}(x_{2}^{2}+...+x_{n}^{2})\leq 1\}$ by   $a^t_1(x_1,x_2,...,x_n)=(2^tx_1,2^{-t}x_2,...,2^{-t}x_n)$. It follows from  Statements~\ref{st-nbh},~\ref{fact} that there exists a homeomorphism $\chi^s_i:N(\tilde{Q}^s_i)\to\mathbb{N}^s$  such that $a^1_1|_{\mathbb N^s}=\chi^s_i a^1 (\chi^s_i)^{-1}|_{\mathbb N^s}$. Denote by  $\chi^s:N(\mathbb L^s)\to \mathbb U\times\mathbb Z_{k^s}$ a homeomorphism such that  $\chi^s|_{N(\tilde{Q}^s_i)}=\chi^s_i$ for any $i\in \{1,...,k^s\}$. Put    $\mathbb{Q}^s=(\mathbb{S}^{n-1}\times \mathbb{R})\cup_{\chi^s}(\mathbb U\times\mathbb Z_{k^s})$. A topological space  $\mathbb{Q}^s$  is a connected oriented $n$-manifold without boundary.  
  
Denote by  $\pi_{s}:(\mathbb{S}^{n-1}\times \mathbb{R})\cup (\mathbb U\times\mathbb Z_{k^s}) \to Q^s$  a natural projection. Put  $\pi_{_{s,1}}=\pi_{_{s}}|_{\mathbb{S}^{n-1}\times \mathbb{R}}$,  $\pi_{_{s,2}}=\pi_{_{s}}|_{\mathbb U\times\mathbb Z_{k^s}}$. 
Define a flow   $\tilde Y^t_s$ on the manifold  $Q^s$ by\\
$\tilde Y^t_{s}(x)=\begin{cases}\pi_{_{s,1}}
(a^t(\pi_{_{s,1}}^{-1}(x))),~
x\in 
\pi_{_{{s,1}}}(\mathbb{S}^{n-1}\times \mathbb{R});
\cr \pi_{_{s,2}}(
a^t_1(\pi_{_{s,2}}^{-1}(x))),~x\in
\pi_{_{s,2}}(\mathbb U\times\{i\}),~i\in\mathbb Z_{k^s}\end{cases}$. 
 
By construction the non-wandering set of the flow    $\tilde Y^t_s$  consists of  $k^s$   equilibria  such that  the flow   $\tilde{Y}^t_s$  is locally  topologically conjugated with the flow $a^t_1$ at  the neighborhood of   each equilibrium. 

{\bf Step 2.} Denote the images of the sets  $\mathbb L^u,~N(\mathbb L^u)$  by means of  the projection $\pi_{_s}$ by the  same  symbols as  their originals. Due to Statements~\ref{st-nbh},~\ref{fact} there exists a homeomorphism   $\chi^u_i:N(\tilde{Q}^u_i)\to \mathbb{N}^u$  such that $a^{-1}_1|_{\mathbb{N}^u}=\chi^u_i \tilde Y^1_s (\chi^u_i)^{-1}$,   $i=1,\dots,k^u$. Denote by  $\chi^u:N(\mathbb L^u)\to\mathbb U\times\mathbb Z_{k^u}$ the  homeomorphism such that
$\chi^u|_{N(\tilde{Q}^u_i)}= \chi^u_i|_{N(\tilde{Q}^u_i)}$  for any $i=1,\dots,k^u$.  Put  $\mathbb{Q}^u=\mathbb{Q}^s\cup_{\chi^u}(\mathbb{U} \times\mathbb Z_{k^u})$. A topological space  $\mathbb{Q}^u$  is a connected oriented $n$-manifold without boundary.
  
Denote by  $\pi_{_{u}}:\mathbb{Q}^s\cup(\mathbb U\times\mathbb Z_{k^u})\to
\mathbb{Q}^u$ the natural projection. Put  $\pi_{_{u,1}}=\pi_{_{u}}|_{\mathbb{Q}^s}$,  $\pi_{_{u,2}}=p_{_{u}}|_{\mathbb{U}\times\mathbb Z_{l^u}}$. 
Define a flow   $\tilde Y^t_u$ on the manifold  $\mathbb{Q}^u$  by\\
$\tilde Y^t_{u}(x)=\begin{cases}\pi_{_{u,1}}
(\tilde Y^t_{s}(\pi_{_{u,1}}^{-1}(x))),~
x\in \pi_{_{{u,1}}}(\mathbb{Q}^s);
\cr \pi_{_{u,2}}(
a^{-t}_1(\pi_{_{u,2}}^{-1}(x))),~x\in
\pi_{_{u,2}}(\mathbb U\times\{i\}),~i\in\mathbb Z_{k^u}\end{cases}$. 

The non-wandering set  $\Omega_{\tilde Y^t_u}$ of the flow  $\tilde Y^t_{u}$  consists of  $k^s$  equilibria  such that the flow   $\tilde{Y}^t_u$ is locally topological  conjugated with the flow  $a^t_1$ in  each of their neighborhoods    and   $k^u$  equilibria such that  the flow  $\tilde{Y}^t_u$ is  locally  topologically conjugated with the flow  $a^{-t}_1$ in each of  their neighborhoods.

{\bf Step 3.} Put $R^s=Q^u\setminus W^s_{\Omega_{\tilde Y^t_u}}$,  denote by  $\rho^s_1,\dots,\rho^s_{n^s}$ connected components of the set  $R^s$ and put     $\hat{\rho}^s_i=\rho^s_i/_{\tilde{Y}^1_u}$. A union of the orbit spaces  $\bigcup\limits_{i=1}^{n^s}\hat{\rho}^s_i$  is obtained from the manifold  $\widehat{V}_{a}$ by a sequence of the surgeries along essential submanifolds of codimension  1. In virtue of  Proposition~\ref{forward} for any  $i\in\{1,...,n^s\}$ the manifold  $\hat{\rho}^s_i$ is homeomorphic to  $\mathbb{S}^{n-1}\times\mathbb{S}^1$,  the manifold  $\rho^s_i$ is homeomorphic to $\mathbb S^{n-1}\times\mathbb R$ and the flow  $\tilde Y^t_u|_{\rho^s_i}$ is topologically conjugated with the flow  $a^t|_{\mathbb R^n\setminus O}$ by means of  a homeomorphism  $\nu^s_i$. Denote by  $\nu^s:R^s\to(\mathbb R^n\setminus \{0\})\times\mathbb Z_{n^s}$ the homeomorphism consisting of  the homeomorphisms  $\nu^s_1,\dots,\nu^s_{n^s}$. Put   ${M^s}=Q^u\cup_{\nu^s}(\mathbb R^n\times\mathbb Z_{n^s})$. Then ${M^s}$ is a connected oriented $n$-manifold without boundary. 
  
Put $\bar{M}^s=Q^u\cup(\mathbb R^n\times\mathbb Z_{n^s})$ and denote by $q_{_{s}}:\bar M^s\to M^s$ the natural projection. Put  $q_{_{s,1}}=q_{_{s}}|_{Q^u}$,  $q_{_{s,2}}=q_{_{s}}|_{\mathbb R^n\times\mathbb Z_{n^s}}$.  Define a flow  $\tilde X^t_s$ on the manifold  ${M^s}$  by \\
$\tilde X^t_{s}(x)=\begin{cases}q_{_{s,1}}
(\tilde Y^t_{u}(q_{_{s,1}}^{-1}(x))),~
x\in q_{_{{s,1}}}(Q^u); \cr 
q_{_{s,2}}(
a^{t}(q_{_{s,2}}^{-1}(x))),~x\in
q_{_{s,2}}(\mathbb R^n\times\{i\}),~i\in\mathbb Z_{n^s}\end{cases}$. 

By construction the non-wandering set of the time one map of the  flow   $\tilde X^t_{s}$ consists of  $k^s$  saddle topologically  hyperbolic fixed points of  index 1,  $k^u$ saddle topologically  hyperbolic fixed points of  index  $(n-1)$ and  $n^s$ sink  topologically hyperbolic fixed points. 
   
{\bf Step 4.} Put $R^u=M^s\setminus W^u_{\Omega_{\tilde X^t_s}}$ and denote by    $\rho^u_1,\dots,\rho^u_{n^u}$  connected components of the set $R^u$. Similar to Step~3 one  can prove that every component $\rho^u_i$ is homeomorphic to  $\mathbb S^{n-1}\times\mathbb R$ and the flow  $\tilde X^t_s|_{\rho^u_i}$ is conjugated with the flow  $a^{-t}|_{\mathbb R^n\setminus \{O\}}$ by a homeomorphism  $\mu^u_i$. Denote by $\mu^u:R^u\to(\mathbb R^n\setminus \{O\})\times\mathbb Z_{n^u}$ a homeomorphism consisting of  the homeomorphisms  $\mu^u_1,\dots,\mu^u_{n^u}$. Put   ${M^u}=M^s\cup_{\mu^u}(\mathbb R^n\times\mathbb Z_{n^u})$. ${M^u}$ is a connected closed oriented     $n$-manifold. 
  
Put $\bar{M}^u=M^s\cup(\mathbb R^n\times\mathbb Z_{n^u})$,  denote by  $q_{_{u}}:\bar M^u\to M^u$ the natural projection, and put $q_{_{u,1}}=q_{_{u}}|_{M^s}$,  $q_{_{u,2}}=q_{_{u}}|_{\mathbb R^n\times\mathbb Z_{n^u}}$. Define a flow  $\tilde X^t_u$ on the manifold ${M^u}$ by \\
$\tilde X^t_{u}(x)=\begin{cases}q_{_{u,1}}
(\tilde X^t_{s}(q_{_{u,1}}^{-1}(x))),~
x\in q_{_{{u,1}}}(M^s);
\cr q_{_{u,2}}(
a^{-t}_0(q_{_{u,2}}^{-1}(x))),~x\in
q_{_{u,2}}(\mathbb R^n\times\{i\}),~i\in\mathbb Z_{n^u}\end{cases}$. 

By construction the non-wandering set of the time one map of the flow   $\tilde X^t_{u}$  consists of  $k^s$  saddle topologically  hyperbolic fixed points of  index 1,  $k^u$ saddle topologically  hyperbolic fixed points of index  $(n-1)$,   $n^s$ sink and $n^u$ source   topologically hyperbolic fixed points. 

{\bf Step 5.} Put $\tilde f=\tilde X^1_u$. By construction $\tilde f$ is a Morse-Smale homeomorphism on the manifold  $M^u$ and its restriction   $\tilde f|_{V_{\tilde f}}$ is topologically conjugated with the diffeomorphism  $f|_{V_f}$ by a homeomorphism mapping the $(n-1)$-dimensional separatrices of the diffeomorphism  $\tilde f$ to the  $(n-1)$-dimensional separatrices of the diffeomorphism  $f$ and preserving their stability.  Due to  Statement~\ref{mat} homeomorphisms $\tilde f$ and  $f$ are topologically conjugated. Hence $M^u=S^n$ and  $X^t=\tilde X^t_u$ is the desired flow.

\end{document}